%% file: git.tex
\documentclass[11pt]{amsart}
\usepackage[dvips]{graphicx,color}
\usepackage{amssymb,amscd,latexsym,epsfig,xypic}
\xyoption{curve}
\usepackage[mathscr]{euscript}
\DeclareFontFamily{OT1}{rsfs}{}
\DeclareFontShape{OT1}{rsfs}{n}{it}{<-> rsfs10}{}
\DeclareMathAlphabet{\curly}{OT1}{rsfs}{n}{it}

\newcommand\C{\mathbb C}

\newcommand\I{\curly I}
\newcommand\g{\mathfrak g}
\renewcommand\k{\mathfrak k}
\newcommand\X{\curly X}
\newcommand\Y{\curly Y}
\newcommand\K{K\"ahler~}
\renewcommand\L{\mathcal L}
\renewcommand\O{\mathcal O}
\newcommand\PP{\mathbb P}
\renewcommand\P{\mathcal P}
\newcommand\Q{\mathbb Q}
\newcommand\R{\mathbb R}
\newcommand\Z{\mathbb Z}
\newcommand\su{\mathfrak su}
\newcommand\tH{\widetilde{\text{Ham}}(X,\omega)}

\newcommand{\rt}[1]{\stackrel{#1\,}{\rightarrow}}
\newcommand{\Rt}[1]{\stackrel{#1\,}{\longrightarrow}}
\newcommand\ot{\leftarrow}
\newcommand\To{\longrightarrow}
\newcommand\into{\hookrightarrow}
\newcommand\Into{\,\ar@{^{(}->}[r]}
\newcommand\acts{\curvearrowright}

\newcommand\dbar{\,\overline{\!\partial}}

\newcommand\ip{\lrcorner\,}

\renewcommand\_{^{}_}
\newcommand\take{\backslash}

\newcommand\rk{\operatorname{rank}}
\newcommand\tr{\operatorname{tr}}

\newcommand\id{\operatorname{id}}
\newcommand\vol{\operatorname{vol}}

\newcommand\Hom{\operatorname{Hom}}

\newcommand\Proj{\operatorname{Proj}\,}
\newcommand\Spec{\operatorname{Spec}\,}
\newcommand\Hilb{\operatorname{Hilb}}
\newcommand\Bl{\operatorname{Bl}}

\renewcommand\sec{\bigskip\noindent\textbf}
\newcommand\beq[1]{\begin{equation}\label{#1}}
\newcommand\eeq{\end{equation}}
\newcommand\beqa{\begin{eqnarray*}}
\newcommand\eeqa{\end{eqnarray*}}

\makeatletter \@addtoreset{equation}{section} \makeatother

\newtheorem{defn}[equation]{Definition}
\newtheorem{thm}[equation]{Theorem}
\newtheorem{lem}[equation]{Lemma}
\newtheorem{cor}[equation]{Corollary}
\newtheorem{conj}[equation]{Conjecture}
\newtheorem{prop}[equation]{Proposition}

\title[Notes on GIT]{Notes on GIT and symplectic reduction for bundles and varieties}
\author{R. P. Thomas}
\date{}

\begin{document}
\maketitle


\section{Introduction}

These notes give an introduction to Geometric Invariant Theory (GIT) and
symplectic reduction, with lots of pictures and simple examples. We describe their applications to moduli of bundles and varieties, leaving the technical
work on the analysis of the Hilbert-Mumford criterion in these situations
to the final sections. We outline their infinite dimensional analogues
(so called Hitchin-Kobayashi correspondences) in gauge theory and in the theory of constant scalar curvature \K (cscK)
and K\"ahler-Einstein (KE) metrics on algebraic varieties. Donaldson's work
on why these should be thought of as the classical limits of the original finite dimensional constructions -- which are then their ``quantisations" -- is explained. The many analogies and connections between the bundle and variety cases are emphasised; in particular the GIT analysis of stability of bundles is shown
to be a special case of the (harder) problem for varieties.

For GIT we work entirely over $\C$ and skip or only sketch many proofs. The interested reader should refer to the excellent \cite{Dl, GIT, Ne} for more details. Throughout this survey we mention many names, but only include certain
key papers in the references -- apologies to those omitted but compiling
a truly comprehensive bibliography would be fraught with danger.
\smallskip

\noindent\textbf{Acknowledgements}. I have learnt GIT from Simon Donaldson
and Frances Kirwan over many years. Large parts of these notes are nothing but an account of Donaldson's point of view on GIT and symplectic reduction for moduli of varieties.
The last sections describe joint work with Julius Ross, and I would like
to thank him, G\'abor Sz\'ekelyhidi and Xiaowei Wang for useful conversations.
Thanks also to Claudio Arezzo, Joel Fine, Daniel Huybrechts, Julien Keller, Dmitri Panov, Michael Singer and Burt Totaro for comments on the manuscript.

\section{A brief review of affine and projective geometry}
This section can be safely skipped by readers with any knowledge of algebraic
geometry. We fix some notation and speedily review some standard theory of complex affine and projective varieties (and schemes). These are much simpler than arbitrary varieties in that they can be described by a single ring. Throughout ``ring" means finitely generated commutative $\C$-algebra -- i.e. a Noetherian ring with a scalar action of $\C$ making it into a $\C$-vector space and commuting with multiplication.

\sec{Affine varieties.}
\emph{Affine varieties} $X$ are just the irreducible components of the zero sets of finite collections of polynomials $p_1,\ldots,p_k$ in some affine space $\C^n$.  They are in one-to-one correspondence with rings with no
zero divisors (i.e. integral domains); in coordinates this is particularly simple:
\begin{eqnarray}
\C[x_1,\ldots,x_n] & \longleftrightarrow & \qquad \C^n, \nonumber \\
\frac{\C[x_1,\ldots,x_n]}{(p_1,\ldots,p_k)} & \longleftrightarrow &
(p_1=0=\ldots=p_k)\subseteq\C^n , \label{corr} \\
\O_X \qquad & \longleftrightarrow & \qquad X. \nonumber
\end{eqnarray}
The arrow $\leftarrow$ replaces a variety by the ring of functions on it (i.e. the functions $\C[x_1,\ldots,x_n]$ on $\C^n$ divided out by the ideal of those that vanish on $X$). The arrow $\to$ recovers $X$ from its ring of functions $\O_X$ by
taking a finite number of generators $x_1,\ldots,x_n$ and a finite (by
the Hilbert basis theorem) number of relations $p_1,\ldots,p_k$ (considered as polynomials in the generators) and setting $X$ to be the affine variety in $\C^n$ cut out by the polynomials $p_i$.

The embedding is equivalent to the choice of generators: each is a map $X\to\C$
so $n$ of them give the map $X\into\C^n$. (Invariantly, we embed $X$ in the
\emph{dual} of the vector space on the generators,
\beq{embed}
X\ni x\mapsto ev_x\in(\C\langle x_1,\ldots,x_n\rangle)^*, \qquad ev_x(f):=f(x),
\eeq 
each point $x$ of $X$ mapping to the linear functional that evaluates functions
at $x$.) The ideal of functions $(p_1,\ldots,p_k)$
vanishing on $X\subset\C^n$ is \emph{prime} (if it contains $fg$ then it contains one of $f$ or $g$) reflecting the fact that the ring $\O_X$ has no zero divisors and that $X$ is irreducible.

So really (\ref{corr}) is a correspondence between
\begin{enumerate}
\item Affine varieties
$X\subseteq\C^n$ with a fixed embedding into $n$-dimensional
affine space,
\item Prime ideals $I\subseteq\C[x_1,\ldots,x_n]$, and
\item Rings without zero divisors plus a choice of $n$ generators.
\end{enumerate}

The coordinate-free approach (which
also shows the above construction is independent of choices) is to note that
the points $x$ of $X$ are in one-to-one correspondence with the maximal ideals
$\I_x\subset\O_X$ of functions vanishing at $x$. So to any ring without zero divisors
$R$ we associate an affine variety $\Spec R$ whose points are maximal ideals in $R$; the coordinate-independent version of (\ref{corr}):
$$
\O_X=R \quad \longleftrightarrow \quad X=\Spec R.
$$
This can be extended to a correspondence between arbitrary rings and affine \emph{schemes}, which are allowed nilpotents in their ring of functions,
corresponding to multiplicities or infinitesimal directions in the scheme.

\sec{Projective varieties.}
Now we just do everything $\C^*$-equivariantly. Recall that a $\C^*$-action
on a vector space $V$ is equivalent to a \emph{grading}, i.e. a splitting into subspaces (the eigen- or weight spaces) $V_k$ parameterised by the integers (the eigenvalues or weights). So we replace rings by graded rings, ideals by homogeneous ideals (those which are the sum of their graded pieces -- i.e. which are $\C^*$-invariant), and get a correspondence between graded rings without zero divisors and projective varieties. This is easiest to
express in coordinates when the ring is generated by its degree one piece:
\begin{eqnarray}
\C[x_0,\ldots,x_n] & \longleftrightarrow & \qquad \PP^n, \nonumber \\
\frac{\C[x_0,\ldots,x_n]}{(p_1,\ldots,p_k)} & \longleftrightarrow &
(p_1=0=\ldots=p_k)\subseteq\PP^n. \label{proj}
\end{eqnarray}
Here $\C[x_0,\ldots,x_n]$ is given the standard grading it inherits from
the scalar $\C^*$-action on $\C^{n+1}$ (i.e. the $x_i$ have weight one), and the $p_i$ are homogeneous polynomials (eigenvectors for the $\C^*$-action).
They cut out an affine variety $\widetilde X=(p_1=0=\ldots=p_k)$ in $\C^{n+1}$, which is $\C^*$-invariant and so a cone, determined by its set of lines through
the origin $X\subseteq\PP^n$.
($X$ and $\PP^n$ are of course the quotients of $\widetilde X\take\{0\}$ and $\C^{n+1}\take\{0\}$ by $\C^*$, but we have yet to develop a theory of quotients (!).)

This describes the arrow $\to$. For $\leftarrow$ we cannot simply take the ring of functions on $X$ since this consists of just the constants; we have to take the ring of functions on the cone $\widetilde X$, which can be described
on $X$ in terms of a line bundle.

Since $X$ is the space of lines in $\widetilde X$, it has a
tautological line bundle $\O_X(-1)=\O_{\PP^n}(-1)|_X$ over it whose fibre
over a point in $X$ is the corresponding line in $\widetilde X\subseteq\C^{n+1}$. The total space of $\O_X(-1)$ therefore has a tautological map to $\widetilde X$ which is an isomorphism away from the zero section $X\subset\O_X(-1)$, which is all contracted down to the origin in $\widetilde X$. In fact the total space of $\O_X(-1)$ is the \emph{blow up} of $\widetilde X$ in the origin.

\medskip \begin{center} \input{blowup.pstex_t} \end{center} \medskip

Linear functions on $\C^{n+1}$ like $x_i$, restricted to $\widetilde X$ and pulled back to the total space of $\O_X(-1)$, give functions which are linear on the fibres, so correspond to sections of the \emph{dual} line bundle $\O_X(1)$. Similarly degree $k$ homogeneous polynomials on $\widetilde X$ define functions on the total space of $\O_X(-1)$ which are of degree $k$ on the fibres, and so give sections of the $k$th tensor power $\O_X(k)$ of the dual of the line bundle $\O_X(-1)$.

So the grading that splits the functions on $\widetilde X$ into homogeneous degree (or $\C^*$-weight spaces) corresponds to sections of different line bundles $\O_X(k)$ on $X$. So $\leftarrow$ takes the direct sum
$\bigoplus_{k\ge0}H^0(\O_X(k))$, considered a graded ring by tensoring sections
$\O(k)\otimes\O(l)\Rt{\simeq}\O(k+l)$. For the line bundle $\O_X(1)$ sufficiently
positive, this ring will be generated in degree one.
It is often called the (homogeneous) coordinate ring of the \emph{polarised}
(=\,endowed with an ample line bundle) variety $(X,\O_X(1))$.

The degree one restriction is for convenience and can be dropped (by working
with varieties in weighted projective spaces), or bypassed by replacing $\O_X(1)$
by $\O_X(p)$, i.e. using the ring $R^{(p)}=\bigoplus_{k\ge0}R_{kp}$; for $p\gg0$ this will be generated by its degree one piece $R_p$.

The choice of generators of the ring is what gives the embedding in projective
space. In fact the sections of any line bundle $L$ over $X$ define a (rational) map
\beq{kod}
X \dashrightarrow \PP(H^0(X,L)^*), \qquad
x \mapsto ev_x, \qquad ev_x(s):=s(x),
\eeq
(compare (\ref{embed})) which in coordinates maps $x$ to $(s_0(x):\ldots:s_n(x))\in\PP^n$, where $s_i$ form a basis for $H^0(L)$. This map is only defined for those $x$ with $ev_x\ne0$, i.e. for which $s(x)$ is not zero for every $s$.

It remains to describe $\to$ in a coordinate-free manner, by noting that
the points of $X$, i.e. lines of $\widetilde X$, are $\C^*$-invariant
subvarieties that are minimal among those which are not the origin in $\widetilde X\subseteq\C^{n+1}$; i.e. homogeneous ideals of the homogeneous coordinate ring that are maximal amongst all homogeneous ideals minus the one corresponding to the origin. So to any graded ring $R$
we associate a projective variety $\Proj R$ whose points are the
homogeneous ideals of $R$ maximal amongst those except the \emph{irrelevant} ideal $R^+:=\bigoplus_{k>0}R_k$. It comes equipped with an ample line bundle, the sections of whose $k$th power gives $R_k$ for $k\gg0$. This gives the
coordinate-free (not quite one-to-one) correspondence (\ref{proj}):
$$
\O_{\widetilde X}=R \quad \longleftrightarrow \quad (X,\O_X(1))=\Proj R.
$$
Replacing $R$ by $R^{(p)}$ leaves the variety $\Proj R$ unaltered but
changes the line bundle from $\O_X(1)$ to $\O_X(p)$.

Similarly arbitrary graded rings with zero divisors (finitely generated, as usual) correspond (not quite one-to-one) to polarised schemes.

\sec{Notation.}
Throughout $G$ will be a connected \textbf{reductive complex linear algebraic group} with Lie algebra $\g$.
Reductive means that all (complex) representations split into sums of irreducibles, but equivalently it is the complexification of a compact real Lie group
$K< G$ with Lie algebra $\k<\g$ such that $\g=\k+i\k$. Therefore representations are also representations of the compact group, which preserve a hermitian inner product (by averaging any inner product on the representation using Haar measure on $K$)
and so split into direct sums of irreducibles by taking orthogonal complements
to invariant subspaces. These splittings are by complex subspaces, so are then also preserved by the complexification $G$.

For the purposes of these notes one can always assume that $K< G$ is
$S^1<\C^*,\ (S^1)^m<(\C^*)^m$ or $SU(m)< SL(m,\C)$.

\section{Geometric Invariant Theory}

GIT is a way of taking quotients in algebraic geometry. This may sound like
a dry and technical subject, but it is beautifully geometric (as we
hope to show) and leads, through its link with symplectic
reduction, to unexpected mathematics (some of which we describe later).

Suppose we are in the following situation, of $G$ acting on a projective
variety $X$ through $SL$ transformations of the projective space.
\beq{G}
\begin{array}{ccc}
G & \acts\quad & X \\
\downarrow && \cap \\
SL(n+1,\C) & \acts\quad & \PP^n.\!\!
\end{array}
\eeq
We would like to form a quotient $X/G$, ideally within the same category
of projective varieties. There are a number of problems with this.

\sec{The topological quotient is not Hausdorff.}
Since $X$ is compact but $G$ is noncompact, a nontrivial $G$-action cannot be proper. There are nonclosed orbits (with lower dimensional orbits in their
closures) so the topological quotient is not Hausdorff.
 
\medskip \hskip 4cm \input{orbits.pstex_t} \medskip

It is clear from the above simple illustration of 3 orbits, all of whose
closures contain the smaller orbit, that we must remove some orbits to get a separated (Hausdorff) quotient. In the above case this would be
the lower dimensional orbit; just as we would expect to remove the origin from $\C^{n+1}$ if we wanted to quotient by the scalar action of $\C^*$ to
get $\PP^n$ (an example to which we shall return).

\sec{Removing smaller orbits does not suffice.}
Another simple example shows that the quotient can still be nonseparated
if we remove all lower dimensional orbits. Consider the action of $\C^*$ on $\C^2$ (or its projective completion $\PP^2\supset\C^2$) by matrices
\vspace{15mm}
\beq{hyp}
\C^*\,\ni\ \lambda\mapsto
\begin{pmatrix}\lambda & 0 \\ 0 & \lambda^{-1}\!\end{pmatrix}
\ \in\,SL(2,\C)
\hskip 5.5cm \vspace{-27mm}
\eeq
\hskip 8cm \epsfig{file=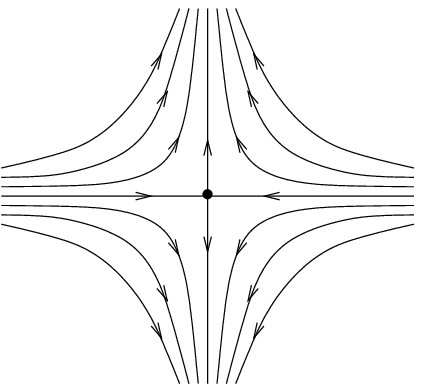}
\medskip

\noindent In this case removing the origin would make the topological quotient the affine line with nonseparated doubled point at the origin: it is clear that the punctured (origin removed) $x$- and $y$-axes are \emph{both} orbits
in the limit of the orbits $\{xy=\alpha\}$ as $\alpha\to0$.

In this simple case it is clear that we would like the quotient to be $\C$, with
$\alpha\in\C\take\{0\}$ representing the orbit $\{xy=\alpha\}$. But what
then should the point $\alpha=0$ of this ``quotient" represent ? There are
three standard solutions to the problem.

\begin{itemize} \item
Kapranov's \textbf{Chow quotient} (and the closely related \textbf{Hilbert
quotient}, of which it is a contraction) parameterises the invariant conics
in this example. For $\alpha\ne0$ these are just the closed orbits of top
dimension, while $\alpha=0$ represents the invariant conic $\{xy=0\}$, the
\emph{union} of all 3 bad orbits. We will not say more about the Chow quotient;
instead see \cite{Hu}, for example (which, unsurprisingly, uses the same
example for illustration).

\item
The \textbf{Geometric Invariant Theory quotient} gives the same for $\alpha\ne0$,
but $\alpha=0$ represents any one of the 3 orbits: GIT identifies all 3 orbits with each other in an equivalence class.

\item
The \textbf{Symplectic reduction} throws away the two nonclosed orbits -- the (punctured) $x$- and $y$-axes -- and keeps the origin. So in this case the quotient represents the closed orbits, including the lower dimensional
one.
\end{itemize}

In this simple case all three quotients give the
same answer $\C$, but in general one expects the Chow quotient to dominate
the GIT and symplectic quotients, which are isomorphic.

In general GIT and the symplectic quotient choose certain ``unstable" orbits to remove to give a separated quotient. GIT also identifies some ``semistable" orbits (those whose closures intersect each other) to compactify the quotient, resulting in a projective variety which we call $X/G$. The symplectic quotient compactifies by taking a distinguished
representative of each semistable equivalence class above -- the intersection of their closures, which turns out to be the unique \emph{closed} semistable orbit in the equivalence class.

\sec{The construction of the GIT quotient.} This is trivial, a formality.
We consider $X$ projective first, since although the affine case is often
even easier, sometimes it is best embedded in the projective case, as we shall see below.

Since we have assumed that $G$ acts through $SL(n+1,\C)$ (rather
than just its quotient $PSL(n+1,\C)$), the action lifts from $X$
to one covering it on $\O_X(-1)$. In other words we don't just act on the
projective space (and $X$ therein) but on the vector space overlying it (and
the cone $\widetilde X$ on $X$ therein). This is called a \emph{linearisation}
of the action. Thus $G$ acts on each $H^0(\O_X(r))$.

Then, just as $(X,\O_X(1))$ is determined by its graded ring of sections
of $\O(r)$ (i.e. the ring of functions on $\widetilde X$),
$$
(X,\O(1))\ \longleftrightarrow\ \bigoplus_r H^0(X,\O(r)),
$$
we simply \emph{construct} $X/G$ (with a line bundle on it) from the ring of \emph{invariant} sections: 
$$
X/G\ \longleftrightarrow\ \bigoplus_r H^0(X,\O(r))^G.
$$
This is sensible, since if there is a good quotient then functions on it
pullback to give $G$-invariant functions on $X$, i.e. functions constant
on the orbits, the fibres of $X\to X/G$. For it to work we need

\begin{lem} \label{fg}
$\bigoplus_r H^0(X,\O(r))^G$ is finitely generated.
\end{lem}

\begin{proof}
Since $R:=\bigoplus_r H^0(X,\O(r))$ is Noetherian, Hilbert's basis theorem
tells us that the ideal $R.\left(\bigoplus_{r>0}H^0(X,\O(r))^G\right)$ generated by $R_+^G:=\bigoplus_{r>0}H^0(X,\O(r))^G$ is generated by a finite number of elements $s_0,\ldots,s_k\in R_+^G$.

Thus any element $s\in H^0(X,\O(r))^G,\ r>0,$ may be written
$s=\sum_{i=0}^kf_is_i$, for some $f_i\in R$ of degree $<r$.
To show that the $s_i$ generate $R^G_+$ as an algebra we must show that the
$f_i$ can be taken to lie in $R^G$.

We now use the fact that $G$ is the complexification of the compact
group $K$. Since $K$ has an invariant metric, we can average over it and use the facts that $s$ and $s_i$ are invariant to give
$$
s=\sum_{i=0}^k\mathop{Av}(f_i)s_i,
$$
where $\mathop{Av}(f_i)$ is the ($K$-invariant) $K$-average of $f_i$. By complex linearity $\mathop{Av}(f_i)$ is also $G$-invariant (for instance,
since $G$ has a polar decomposition $G=K\exp(i\k)$). The $\mathop{Av}(f_i)$ are also of degree $<r$, and so we may assume, by an induction on $r$, that we have already shown that they are generated by the $s_i$ in $R_+^G$. Thus $s$ is also.
\end{proof}

Thus we simply \emph{define} $X/G$ to be $\Proj\bigoplus_r H^0(X,\O(r))^G$.
If $X$ is a variety (rather than a scheme) then so is $X/G$, as its graded
ring sits inside that of $X$ and so has no zero divisors. Unfortunately
this is not all there is to GIT, however. We have to work out what $X/G$
is, which orbits points of $X/G$ represent, and so on, which we tackle in the next section. Another important question, that we barely touch on, is how the quotient $X/G$ changes with the linearisation.
For some linearisations the quotient is empty, but if under a change of
linearisation the moduli space remains the same dimension then it undergoes only a small birational transformation, a type of \emph{flip} \cite{DH, Th}.

\sec{The affine case.} The affine case is even easier; if $G$ acts on $\Spec
R$ we can form $\Spec(R^G)$ as a putative quotient. For instance in our example
(\ref{hyp}) the ring of invariants
$$
R^G=\C[x,y]^{\C^*}=\C[xy]
$$
is generated by $xy$, so that the quotient is $\Spec\C[xy]\cong\C$, as anticipated.
The function $xy$ does not distinguish between any of the ``bad" orbits (the punctured
$x$- and $y$- axes, and the origin), lumping them all in an equivalence class
of orbits which get identified in the quotient. \smallskip

In other cases this does not work so well; for instance under the scalar action of $\C^*$ on $\C^{n+1}$ the only invariant polynomials in
$\C[x_0,\ldots,x_n]$
are the constants and this recipe for the quotient gives a single point. In the language of the next section, this is because there are no stable points in this example,
and all semistable orbits' closures intersect (or equivalently, there is a unique polystable point, the origin). More generally in any affine case
all points are always at least semistable (as the constants are always $G$-invariant
functions) and so no orbits gets thrown away in making the quotient (though
many may get identified with each other -- those whose closures intersect
which therefore cannot be separated by invariant functions). But for the scalar action of $\C^*$ on $\C^{n+1}$ we clearly need to remove at least the origin to get a sensible quotient.

So we should change the linearisation,
from the trivial linearisation to a nontrivial one, to get a bigger quotient.

\sec{Example: $\PP^n$ from GIT.}
That is, we consider the trivial line bundle on $\C^{n+1}$ but with a nontrivial linearisation, by composing the $\C^*$-action on $\C^{n+1}$ by a character $\lambda\mapsto\lambda^{-p}$
of $\C^*$ acting on the fibres of the trivial line bundle over $\C^{n+1}$. The invariant sections of
this no longer form a ring; we have to take the direct sum of spaces
of sections of \emph{all powers} of this linearisation, just as in the projective
case, and take Proj of the invariants of the resulting graded ring.

If $p<0$ then there are no invariant sections and the quotient is empty.
We have seen that for $p=0$ the quotient is a single point. For $p>0$ the
invariant sections of the $k$th power of the linearisation are the homogeneous
polynomials on $\C^n$ of degree $kp$. So for $p=1$ we get the quotient
\beq{Pn}
\C^{n+1}/\C^*=\Proj\bigoplus_{k\ge0}\big(\C[x_0,\ldots,x_n]_k\big)=
\Proj\C[x_0,\ldots,x_n]=\PP^n.
\eeq
For $p\ge1$ we get the same geometric quotient but with the line bundle $\O(p)$ on it instead of $\O(1)$.

Another way to derive this is to
embed $\C^{n+1}$ in $\PP^{n+1}$ as $x_{n+1}=1$, act by $\C^*$
on the latter by diag$\big(\lambda,\ldots,\lambda,\lambda^{-(n+1)}\big)
\in SL(n+2,\C)$, and do projective GIT. This gives, on restriction to $\C^{n+1}\subset\PP^{n+1}$, the $p=n+1$ linearisation
above. The invariant sections of $\O((n+2)k)$ are of the form $x_{n+1}^k.f$, where $f$ is a homogeneous polynomial of degree $(n+1)k$ in $x_1,\ldots,x_n$.
Therefore the quotient is
$$
\Proj\bigoplus_{k\ge0}\big(\C[x_1,\ldots,x_n]_{(n+1)k}\big)=
(\PP^n,\O(n+1)).
$$
In the language of the next section this is because the complement of
$\C^{n+1}\subset\PP^{n+1}$, and the origin $\{0\}\in\C^{n+1}\subset\PP^{n+1}$, are unstable (either by noting that all of the nonconstant
invariant polynomials above vanish on them, or by an easy exercise in using
the Hilbert-Mumford criterion below -- these loci are fixed points, but with a nontrivial action on the line above them). So these are removed and the projective quotient reduces to the affine case.

\sec{What are the points of a GIT quotient ?}
By its very definition (and Lemma \ref{fg}), for $r\gg0$, $X/G$ is just the image of $X$ under the linear system $H^0(\O_X(r))^G$.
That is, consider the Kodaira ``embedding" of $X/G$ (\ref{kod}),
\begin{eqnarray} \label{Kod}
\qquad X &\dashrightarrow& \PP((H^0(X,\O(r))^G)^*), \\
x &\mapsto& \quad ev_x \hspace{4cm} (ev_x(s):=s(x)), \nonumber
\end{eqnarray}
that in coordinates takes $x$ to $(s_0(x):\ldots:s_k(x))\in\PP^k$
(where the $s_i$ form a basis for $H^0(X,\O(r))^G$).

This is only a rational map, since it is only defined on points for which
$ev_x\not\equiv0$ (equivalently the $s_i(x)$ are not all zero). That is,
it is defined on the \emph{semistable points} of $X$:

\begin{defn}
$x\in X$ is \textbf{semistable} if and only if there exists $s\in H^0(X,\O(r))^G$
with $r>0$ such that $s(x)\ne0$.
\end{defn}

Points which are not semistable we call (controversially) \textbf{\emph{unstable}}.

So semistable points are those that the $G$-invariant functions ``see". The map (\ref{Kod}) is well defined on the (Zariski open, though possibly empty) locus $X^{ss}\subseteq X$ of semistable points, and it is clearly constant on $G$-orbits, i.e. it factors through the set-theoretic quotient $X^{ss}/G$. But it may contract more than just $G$-orbits, so we need another definition.

\begin{defn}
A semistable point $x$ is \textbf{stable} if and only if $\bigoplus_rH^0(X,\O(r))^G$ separates orbits near $x$ and the stabiliser of $x$ is finite.
\end{defn}

By ``separates orbits near $x$" we mean the following. Since $x$ is semistable
there exists an $s\in H^0(\O_X(r))^G$ such that $s(x)\ne0$. So now we work
on the open set $U\subset X$ on which $s\ne0$ and use $s$ to trivialise $\O_U(r)$
(i.e. divide all sections of $\O_U(r)$ by $s$ to consider them as functions).
Then we ask that in $U$ any orbit can be distinguished from $G.x$ by $H^0(X,\O(r))^G$.
That is, there is an element of $H^0(X,\O(r))^G$ which takes different values
on the two orbits, and this should also be true infinitesimally: given a
vector $v\in T_xX\take T_x(G.x)$, there is an element of $H^0(X,\O(r))^G$ whose derivative down $v$ is nonzero.

So we have a (surjective) map $X^{ss}\to X/G$ under which the line bundle
on $X/G$ (that arises from its Proj construction) pulls back to $\O_{X^{ss}}(1)$. This map has good geometric
properties over the locus of stable points $X^s\subseteq X^{ss}\subseteq X$ (it only contracts single orbits, for instance; more properly
it is a \emph{geometric quotient} in the sense of \cite[Definition 0.6]{GIT}). The definitions of stable and semistable are the algebraic ones designed to make this true, but now we can relate them more to geometry.

\sec{Topological characterisation of (semi)stability.} If we work upstairs in the
vector space $\C^{n+1}\supset\widetilde X$ (or equivalently in the total
space of $\O_X(-1)$) instead of in the projective space $\PP^n\supset X$,
we can get a beautiful topological characterisation of (semi)stability.
Given our topological discussion about nonclosed orbits at the start of these notes, it is what one might guess, and the best one could possibly hope for.
So for $x\in X$, pick $\tilde x\in\O_X(-1)$ covering it.

\begin{thm} \label{Mum} \ \\
$x$ is semistable $\iff 0\not\in\,\overline{\!G.\tilde x\,}$. \\
$x$ is stable $\iff G.\tilde x$ is closed in $\C^{n+1}$ and $\tilde x$ has
finite stabiliser.
\end{thm}

\noindent When $G.\tilde x$ is closed, but not necessarily of full dimension, we call $x$ \textbf{polystable}. (This is called
``Kempf-stable" in \cite{Dl}, ``weakly stable" in \cite{Ti2} and plain ``stable" in \cite{Do4}. Originally \cite{GIT} the terms for stable and polystable
were ``properly stable" and ``stable" respectively.)

In one direction the theorem is clear. $G$-invariant homogeneous functions of degree $r>0$ on $\widetilde X$ are constant on orbits and so also their closures. So if the closure of the orbit of $\tilde x$ contains the origin then every such function is zero on $\widetilde x$ and $x$ cannot be semistable. Similarly if the invariant functions separate orbits around the orbit of $\tilde x$ then it is the zero locus of a collection of invariant functions and so closed.

One can make the criterion (\ref{Mum}) much simpler to calculate with, by
considering one parameter subgroups instead of all of $G$.

\sec{The Hilbert-Mumford criterion.} The key result is that $x$ is (semi)stable for $G$ if and only if it is (semi)stable \emph{for all one parameter subgroups (1-PSs) $\C^*< G$}. We will outline a proof of this remarkable result
once we have done some symplectic geometry.

So we may apply Theorem \ref{Mum} to each of these 1-PS orbits, and determining
the closedness of these one dimensional orbits is much easier by using their asymptotics. Setting
$x_0=\lim_{\lambda\to0}\lambda.x$, this is a fixed point of the $\C^*$-action, so $\C^*$ acts on the line $\O_{x_0}(-1)$ in $\C^{n+1}$ that $x_0\in\PP^n$ represents. Letting $\rho(x)\in\Z$ denote the weight of this action (i.e.
$\C^*\ni\lambda$ acts on $\O_{x_0}(-1)$ as $\lambda^{\rho(x)}$) we find the
following, the Hilbert-Mumford criterion.

\begin{thm}\ \vspace{-4mm} \\
\begin{itemize}
\item If $\rho(x)<0$ for all 1-PS then $x$ is stable,
\item If $\rho(x)\le0$ for all 1-PS then $x$ is semistable,
\item If $\rho(x)>0$ for a 1-PS then $x$ is unstable.
\end{itemize}
\end{thm}

The proof is the picture below; the 1-PS orbit is closed if and only if it is asymptotic
to a negative weight $\C^*$-action on the limiting line at both $\lambda\to0$
and $\lambda\to\infty$. But we can restrict to the former since the
latter arises from the inverse 1-PS.

\bigskip \input{smallHM.pstex_t} \bigskip

So we ``just" have to compute the weight
$\rho(x)$ for \emph{all} $\C^*< SL(n+1,\C)$; $x$ is stable for $G$
if and only if $\rho(x)$ is always $<0$.

To sum up; the 1-PS orbits are not (in general) closed in the projective
space, but they may be upstairs in the vector space. To decide if $x$ is
stable we first take the limit $x_0$ as we move through isomorphic objects $\lambda.x$; this limit is \emph{not} (in general) isomorphic to $x$ under the $\C^*$-action; it is only in the \emph{closure} of the $\C^*$-orbit. This point $x_0$ represents a line $\O_{x_0}(-1)$ in the vector space, on which $\C^*$ acts. If this weight is negative then $x$ (\emph{not} $x_0$!)
is stable for this 1-PS; if this is true for all 1-PSs then $x$ is stable
for $G$.

\sec{Fundamental example: points in $\PP^1$}.
The standard example, from for instance \cite{GIT}, is to consider configurations
of $n$ (unordered) points in $\PP^1$ up to the symmetries of $\PP^1$. (This is of course a 0-dimensional
algebraic variety, and so the easiest example of the stability of varieties
that we shall study later.) In fact we allow multiplicities, i.e. we take
length-$n$ 0-dimensional subschemes of $\PP^1$ -- $S^n\PP^1$ -- modulo $SL(2,\C)$.

To linearise the action we note that specifying any such $n$ points is the
same as specifying a degree $n$ homogeneous polynomial on $\PP^1$, unique up to scale, by taking the roots of the polynomial. That is, $S^n\PP^1$ is
the projectivisation of $H^0(\O_{\PP^1}(n))$, so giving us a natural linearisation
of the problem; we use the induced $SL(2,\C)$-action on $H^0(\O_{\PP^1}(n))
\cong S^n(\C^2)^*$.

We find that the configuration is stable unless it has a very singular point.

\begin{thm} \label{P1}
A length-$n$ subscheme of $\PP^1$ is
\begin{itemize}
\item semistable if and only if each multiplicity $\le n/2$,
\item stable if and only if each multiplicity $<n/2$.
\end{itemize}
\end{thm}

\begin{proof}
1. Diagonalise a given $\C^*< SL(2,\C):$
\beq{diag}
\begin{pmatrix}\lambda^k & 0 \\ 0 & \lambda^{-k}\!\end{pmatrix}
\ \mathrm{in}\ [x:y]\ \mathrm{coordinates\
on}\ \PP^1.\ \,(k\ge0.)
\eeq
In these coordinates write our degree $n$ homogeneous polynomial (whose roots
give the $n$ points) as $f=\sum_{i=0}^na_ix^iy^{n-i}$. As $\lambda\to0$ the
first half (precisely the first $\lceil n/2-1\rceil$) of these monomials tend to infinity (as there are more $y$s than $x$s in the monomial).

Thus $\lambda.f$ tends to $\infty$ and the orbit is closed about $\lambda\to0$ unless $a_i=0$ for $i\le n/2$.
That is, it is closed so long as $f$ does not vanish to order $\ge n/2$ at $x=0$.

Repeating over all 1-PS changes the coordinates $[x:y]$, so $f$ is stable if and only if it does not vanish to order $\ge n/2$ at any point.
\end{proof}

Alternatively, we can use the Hilbert-Mumford criterion in terms of the
weight on the limiting line.

\begin{proof}
2. Up to rescaling, under the action (\ref{diag}), $\lambda.f\to
f_0=a_jx^jy^{n-j}$, where $j$ is smallest such that $a_j\ne0$.

The weight of (\ref{diag}) on $\C.f_0$ is $k(j-(n-j))=k(2j-n)$.
So $f$ is stable if and only if $k(2j-n)<0\iff j<n/2\iff\text{ord}_{\,x=0}(f)<n/2$ for all 1-PS (and so all points $x=0$) as before. Semistability is similar.
\end{proof}

\medskip \begin{center} \epsfig{file=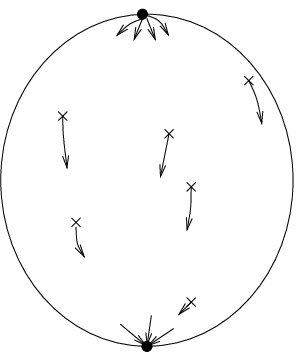} \end{center}\medskip

Geometrically what is happening is that the 1-PS moves almost all points to the ``attractive" fixed point at $x=0$ (weight $-k$), and this is the
generic, stable, situation. Only those points stuck at ``repulsive" fixed point $y=0$ contribute positively ($+k$) to the weight. So the total weight is negative unless more than half of the points are at the repulsive fixed
point. \medskip

More generally if we consider hypersurfaces $D\subset\PP^n$ modulo the action
of $SL(n+1,\C)$, then the existence of the \emph{discriminant} of the defining equation (which is automatically an invariant polynomial) means that
\beq{smooth}
D \text{ smooth }\Rightarrow D \text{ semistable},
\eeq
since its discriminant does not vanish.

\sec{Another topological characterisation of stability.} 
Instead of working in $\widetilde X$, we can give a topological characterisation of stability (not semistability) \emph{downstairs} in $X$ which, in the
language of the first section, says that
we have removed enough orbits to get a Hausdorff quotient. 

Namely, polystability of $x$ is equivalent to the orbit of $x$  being \emph{closed in the locus of semistable points}.

For stability, one direction is clear; by the original definition of stability of $x$, the \emph{closed}
locus where the invariant sections $\bigoplus_kH^0(L^k)^G$ take the value that they take on $G.x$ is precisely $G.x$.

For the converse we note that given any 1-PS orbit of a stable point $x$,
its limit point $x_0$ is \emph{unstable}, by the Hilbert-Mumford criterion
applied to the inverse 1-PS (under which $x_0$ is fixed and $\O_{x_0}(-1)$ is acted on with positive weight). Thus once $x_0$ is removed, i.e. in the locus of semistable points, the 1-PS orbit of $x$ is closed (at the $\lambda\to0$
end; then one can consider the inverse 1-PS). Then one has to show
that this enough to show that the whole of $G.x$ is closed in $X^{ss}$.

The polystable case follows by showing that for any semistable point $x$,
$$
\bigcap\big\{\overline{G.y}\ \colon\ \overline{G.y}\ \cap\ 
\overline{G.x}\ne\emptyset\big\}
$$
is a \emph{single orbit}, which is therefore closed. This is the unique polystable
representative of the semistable equivalence class of orbits which are identified together in the GIT quotient. This can be proved using the Kempf-Ness theorem
in a later Section.

\sec{Another form of the Hilbert-Mumford criterion.} The Hilbert-Mumford
criterion can be recast in terms of $\tilde x$ only (with no mention of $x_0$)
in a form that can be useful in calculations.

In our usual set-up of $G$
acting on $X\subseteq\PP^n$ linearised over $\O_X(-1)$, any 1-PS therefore
gives a 1-PS of $SL(n+1,\C)$ (in fact $GL(n+1,\C)$, but it is an easy exercise
to see that if it does not lie in $SL(n+1,\C)$ then all points are unstable).

Diagonalising we can write the action as $\lambda\mapsto$\,diag\,$(\lambda^{\rho_0},
\ldots,\lambda^{\rho_n})$, in which basis $\tilde x=(x_0,\ldots,x_n)$. Then
it is clear that the orbit $\{\lambda.\tilde x\}$ tends to $\infty$ as $\lambda\to0$
if and only if there is an $i$ such that
\beq{-ve}
x_i\ne0 \quad\text{and}\quad \rho_i<0.
\eeq
This is true for all 1-PSs if and only if $\tilde x$ is stable. The semi-
and poly- stable cases are left as an exercise for the reader.

So the vector space $\C^{n+1}=W^+\oplus W^0\oplus W^-$ can be split into a sum of positive and
negative weight spaces (with the sum of all weights, with multiplicities,
zero) for the $\C^*< SL(n+1,\C)$. The condition that $\tilde x$ be
stable with respect to this 1-PS and its inverse is then that its components in both of $W^+$ and $W^-$ be nonzero; i.e. that it have components of both positive and negative
weight. Hence \emph{generic} points are stable if the group
acts effectively through $SL(n+1,\C)$.

\section{Symplectic reduction}

In this section we use the compact subgroup $K< G$ to enlarge (\ref{G})
to
\beq{K}
\begin{array}{ccccc}
K & < & G & \acts\quad & X \\
\downarrow && \downarrow && \cap \\
SU(n+1) &  < & SL(n+1,\C) & \acts\quad & \PP^n.\!\!
\end{array}
\eeq
So $K$ acts on $\PP^N$, preserving the complex structure $J$ (as $SL(n+1,\C)$
does) but also the Fubini-Study metric $g$. Thus $K$ preserves the symplectic
form $\omega=g(\ \cdot\ ,J\ \cdot\ )$ and acts through \emph{symplectomorphisms}
in Aut$(X,\omega)$.

\sec{Hamiltonian automorphisms.}
The Lie algebra of the symplectomorphism group consists of vector fields
$Y$ on $X$ which preserve the symplectic form, i.e.
$$
\mathcal L\_Y\omega=d(Y\ip\omega)=0.
$$
Since contraction with $\omega$ is an isomorphism $TX\Rt{\simeq}T^*X$, the Lie algebra is isomorphic to the space of \emph{closed 1-forms} $Z^1(X)$.

The subspace of \emph{exact 1-forms} $dC^\infty(X,\R)$ (which is all of them
for simply connected spaces like $\PP^N$) generate \textbf{hamiltonian}
automorphisms -- those which can be connected to the identity through a path
of automorphisms whose flux homomorphism is zero (i.e. integration of $\omega$
over the cylinder traced out in this path by any loop in $X$ is zero).

Since $d$ has kernel the constants, the Lie algebra of the hamiltonian automorphisms
is $C^\infty(X,\R)/\R$; the function $h$ generating the vector field
$$
X_h \quad\text{such that}\quad X_h\ip\omega=dh.
$$
That is, using the metric, $X_h$ is the symplectic gradient $J\nabla h$ of the \emph{hamiltonian} $h$.

The Lie bracket of hamiltonian vector fields works out to be the \textbf{Poisson bracket} $\{f,g\}$ on functions, given by pairing $df$ and $dg$ using the
(inverse of) the symplectic form and dividing out by constants. (Equivalently, this is the class of $X_f(g)=-X_g(f)$ in $C^\infty(X,\R)/\R$.)

This bracket clearly lifts to $C^\infty(X,\R)$, and can be checked to satisfy
the Jacobi identity there too. The constants are central -- they Poisson commute with
all of $C^\infty(X,\R)$ -- so we get a \emph{central extension} of the Lie algebra of hamiltonian automorphisms:
\beq{ext}
0\to\R\to C^\infty(X,\R)\to C^\infty(X,\R)/\R\to0.
\eeq

One might ask what $C^\infty(X,\R)$ is the Lie algebra of, or what it acts on infinitesimally. Since $\R$ is the Lie algebra of isometries of the line $\C$, one might consider isometries of a line bundle $L\to X$ covering hamiltonian automorphisms on $X$.

This indeed can be made to work if $\omega$ is integral, i.e. its cohomology
class lies in $H^2(X,\Z)/$torsion $\le H^2(X,\R\,)$ (as in our projective
case, for instance). Then $2\pi i\omega$ is the curvature of a hermitian
line bundle $L$ with unitary connection, and we let $\tH$ be the isometries
of $L$ preserving its connection; these then cover hamiltonian automorphisms
on $X$. Infinitesimally $h\in C^\infty(X,\R)$ acts through vector fields on $L$ given by
\beq{???}
\tilde{X_h}+ih.
\eeq
Here $\tilde{X_h}$ is the horizontal lift of the hamiltonian vector field
$X_h$, and $ih$ is the multiplication operator taking element of the line
$L$ to a perpendicular element in its tangent space (using the natural isomorphism
between a line $L$ and its tangent space). This action defines a homomorphism of Lie algebras (the Poisson bracket on
$C^\infty(X,\R)$ maps to the Lie bracket on vector fields on the total space
of $L$), and the constants $\R$ act as the Lie algebra
of global constant rotations $\{e^{i\theta}\}=U(1)$ of the fibres of $L$,
yielding the exact sequence (\ref{ext}).

(This is often called \emph{prequantisation}, giving a representation of
the hamiltonian diffeomorphisms on the projectivisation of the Hilbert space of $L^2$-sections of $L$. This is considered too big a Hilbert space to be
the set of quantum mechanical wave functions, and geometric quantisation
attempts to replace it with holomorphic sections; about which more later.)

\sec{Moment maps and linearisations.}
Since $K$ acts through symplectomorphisms of $\PP^n$, which is simply connected
(so that $Z^1(X)=dC^\infty(X)$), we get a Lie algebra homomorphism
$$
\k=\operatorname{Lie}(K)\to C^\infty(\PP^n,\R)/\R\to C^\infty(X,\R)/\R, \qquad v\mapsto[m_v],
$$
where any $v\in\k$ generates a hamiltonian vector field $X_v$ on $X$ such that $X_v\,\ip\omega=dm_v$, for some function $m_v$ unique up to a constant.

We would like to choose these constants consistently, i.e. choose a lift
\beq{lift}
\diagram
& C^\infty(X,\R)\dto &&& m_v\! \ar@{|->}[d] \\
\k \urto\rto & C^\infty(X,\R)/\R
&& v \ar@{|->}[ur]\ar@{|->}[r] & [m_v]\cong X_v=dm_v\hspace{-23mm}
\enddiagram \hspace{2cm}
\eeq
which is a \emph{homomorphism of Lie algebras}. One such always exists, since
(\ref{ext}) is split by the Poisson subalgebra $C^\infty(X,\R)_0\cong C^\infty(X,\R)/\R$
of functions of integral zero; i.e. we can choose each $m_v$ to have integral
zero. But we want to consider arbitrary lifts since from (\ref{???}) we know that each is equivalent to a lift
of $\k$ to isometries of $L$ preserving the connection and covering its hamiltonian action downstairs. This is the infinitesimal version of a linearisation,
assigning to $v\in\k$ the vector field
\beq{imv}
\tilde X_v+im_v
\eeq
on the total space of $L$. It may or may not integrate up to an action of $K$ on $(X,L)$ covering that on $X$.

Said differently, we want to put together all of the hamiltonians $m_v$ to
give a \textbf{moment map}
\beq{mm}
m\,\colon\,X\to\k^*,
\eeq
such that $\langle m(x),v\rangle=m_v(x)$ for all $v\in\k$. $m$ is just
a collection of $\dim K$ hamiltonians $m_v$, written invariantly. Then
our lifting condition (\ref{lift}) becomes the condition that the undetermined constants in $m_v$ be chosen such that (\ref{mm}) is \emph{$K$-equivariant} (using the coadjoint action on the right hand side). Thus a moment map is unique up to the addition of a central element of $\k^*$.

Yet another way of saying the same thing is that the derivative of the $K$-action
maps $\k$ to $TX$, so by contraction with the symplectic form $\ip\omega\colon
TX\cong T^*X$ is a section of $\k^*\otimes T^*X$. It is closed and $K$-invariant, so we ask for it to be invariantly exact, i.e. $d$ of a $K$-invariant section $\mu$ of $C^\infty(X,\k^*)$.

(The name comes from the case of a cotangent bundle $X=T^*M$ with its canonical symplectic form and action induced from an action
of $K$ on $M$. Then the moment map really gives the momentum of the image $X_v\in TM$ of $v\in\k$: $m_v(p,q)=\langle p,X_v\rangle$ at a point $q\in M$ and $p\in T^*_qM$. Hence for translations we get the usual linear
momentum, and for rotations angular momentum.)

In the projective case that we have been considering, a natural $m$ exists because we picked a linearisation. $SU(n+1)\acts(\PP^n,\O(1))$ has a canonical moment map given by
\beq{mmpn}
\tilde x\mapsto\frac{i\big(\langle\ \cdot\ ,\tilde x\rangle\otimes\tilde x\big)_0}{||\tilde x||^2}\in\su(n+1)^*\cong\su(n+1),
\eeq
where $(\ )_0$ denotes the trace-free part of an endomorphism. Restricting to $X$ and projecting to $\k^*$ by the adjoint of the map $\k\to\su(n+1)$ gives a moment map for the $K$-action on $X$.

\sec{The Kempf-Ness theorem.} The key to the link between symplectic
geometry and GIT is the following calculation. Suppose $(X,L=\O_X(1))$ is a polarised
variety with a hermitian metric on $L$ inducing a connection with curvature $2\pi i\omega$. Lift $x$ to any $\tilde x\in\O_x(-1)=L^{-1}_x$ and consider the \textbf{norm functional} $||\tilde x||$. (If $X$ is embedded in $\PP(H^0(L)^*)$
then one way to get a metric
on $\O(-1)$ is to induce it from one on $H^0(L)^*$ upstairs; then $||\tilde x||$ is just the usual norm in the vector space that $\widetilde X$ lives
in.) As we move down a 1-PS orbit $\{\lambda.\tilde x\colon\lambda\in\C^*\}$ in the direction of $v\in\k$ we see how $\log||\tilde x||$ varies;
for $\lambda\in U(1)<\C^*$ (which preserves the metric) not at all, but for $\lambda$ in the complexified, radial direction $\lambda\in(0,\infty)
<\C^*$ we get
\beq{dd}
m_v=\left.\frac d{d\lambda}\right|_{\lambda=1}\log||\lambda\tilde x||
\_{\lambda\in(0,\infty)}\ .
\eeq
That is, $X_v(\log||\lambda\tilde x||)=0$, but
\beq{imv2}
(JX_v)(\log||\lambda\tilde x||)=X_{iv}(\log||\lambda\tilde x||)=m_v.
\eeq
(This is just an unravelling of (\ref{imv}).
For instance if $x$ is a fixed point, then $\C^*$ acts on the line $\langle
\tilde x\rangle$ with a weight $\rho$, and
\beq{rho}
m_v=\rho,
\eeq
which is therefore an integer.)

Moreover, $\log||\lambda\tilde x||$ is \emph{convex} on $\C^*/U(1)\cong(0,\infty)$, as its second derivative is positive:
$$
X_{iv}m_v=dm_v(JX_v)=\omega(X_v,JX_v)=||X_v||^2.
$$
It follows that the orbit tends to infinity at both ends, i.e. is closed,
if and only if it contains a critical point (i.e. absolute minimum) of $\log
||\lambda\tilde x||$.

\smallskip \begin{center} \input{smallorb.pstex_t} \end{center} \smallskip

So a 1-PS orbit is polystable if and only if it contains
a zero of the corresponding hamiltonian. That zero is then unique, up to
the action of $U(1)$. This is the Kempf-Ness theorem for
$\C^*$-actions.

Next we would like to consider a full $G$ orbit, and find a zero of all the hamiltonians simultaneously, i.e. a zero of $m$.
Pick $v_i$ to form a basis for the Lie algebra of a maximal torus in $K$
such that each generates a 1-PS. If an orbit is polystable then each 1-PS
orbit is closed, so by the above there is a point with $m_{v_1}=0$ in the
first 1-PS. Now we move down the second 1-PS orbit of this point to a
point with $m_{v_2}=0$ \emph{and} $m_{v_1}=0$ since the two 1-PSs commute
(i.e. $\{m_{v_1},m_{v_2}\}=0$). Inductively we find a point with $m_v=0$
for all $v$ in the Lie algebra of the torus, and so for all $v$ conjugate
to such (i.e. all $v$) by equivariance of the moment map. Thus the orbit contains a point with $m=0$. Moreover, by the convexity of $\log||\tilde
x||$ on $G/K$, the zero is in fact unique up to the action of $K$.

(Alternatively, we could have proved this without using the Hilbert-Mumford
criterion by noting that $\log||g.\tilde x||$ is convex on the whole of $G/K$,
instead of each $\C^*/U(1)$, so an orbit is closed if and only if this functional
has a minimum, at which point $m=0$ by (\ref{dd}).)

\begin{thm}\emph{[Kempf-Ness]}
A $G$-orbit contains a zero of the moment map if and only if it is polystable.
It is unique up to the action of $K$.

A $G$-orbit is semistable if and only if its closure contains a zero of the
moment map; this zero is in the unique polystable orbit
in the closure of the original orbit.

In particular, as sets,
$$
\frac XG\ \cong\ \frac{m^{-1}(0)}K\,=:X/\kern-.7ex/K.
$$
\end{thm}

$X/\kern-.7ex/K:=m^{-1}(0)/K$ is called the \textbf{symplectic reduction} of $X$, invented by Marsden-Weinstein and Meyer. \medskip

\begin{center}
\input{slice4.pstex_t}
\end{center} \smallskip

So on the locus of stable points $m^{-1}(0)$ provides a ($K$-equivariant) slice to the $i\k<\g=\k+i\k$ part of orbit; since this is topologically trivial ($G$ retracts onto $K$) it makes topological sense that one could take a slice instead of a quotient. This leaves only the $K$-action to divide by to get the GIT quotient.

The Kempf-Ness theorem is a nonlinear generalisation of the isomorphism
$V/W\cong W^\perp$ for vector spaces $W\le V$. It works due to convexity,
giving a unique distinguished $K$-orbit of points of least norm in each polystable $G$-orbit upstairs in $\widetilde X$.

When 0 is a regular value of $m$ (which implies that $m^{-1}(0)$ is smooth and the $\k$-action on it is injective, so the $K$-action has finite stabilisers
and the quotient is a smooth orbifold at worst) then the restriction of $\omega$ to $m^{-1}(0)$ is degenerate precisely along
the $K$-orbits, and so descends to a symplectic form on the quotient. This
is in fact compatible with the complex (algebraic) structure on the GIT quotient,
giving a \K form representing the first Chern class of the polarisation that $X/G$ inherits from its Proj construction.

\sec{Example.} $U(1)<\C^*$ acts on $\C^{n+1}$ with moment map $m=
|\underline z|^2-a$ for any constant $a\in\R$. For $a>0$ this gives
$$
\frac{\C^{n+1}\take\{0\}}{\C^*}\ \cong\ 
\frac{S^{2n+1}=\{\underline z\,:\ |\underline z|^2=a\}}{U(1)}\ \cong\ \PP^n.
$$
$S^{2n+1}=m^{-1}(0)$ is a slice to the $(0,\infty)$-action, leaving the $U(1)$-action
to divide by. The resulting \K form on $\PP^n$ varies with the level
$a$.

For $a=0$ we get just a single point, while for $a<0$ we get the empty set
-- as we showed already using GIT for different polarisations (\ref{Pn}), where $p$ played the role of $a$ (but took integer values so that the lifted action of $\k$ descended to an action of $K=U(1)$ on the trivial line bundle
over $\C^{n+1}$).

\sec{Example: $n$ points in $\PP^1$ again.} (Kirwan \cite{Ki}) The moment
map
$$
SL(2,\C)\supset SU(2)\acts\PP^1\stackrel{m\,}\To\su(2)^*
$$
is just the inclusion of the unit sphere $S^2\subset\R^3$.

Adding gives, for $n$ points, the moment map $m=\sum_{i=1}^nm_i$:
\beq{mm1}
S^n\PP^1\To\R^3,
\eeq
the sum of the $n$ points in $\R^3$, i.e. ($n$ times) their centre of mass.

So $m^{-1}(0)$ is the set of \textbf{balanced configurations} of points with
centre of mass $0\in\R^3$.

Since by Kempf-Ness polystability is equivalent to the existence of an $SL(2,\C)$ transformation of $\PP^1$ that balances the points, Theorem \ref{P1} yields

\begin{thm} A configuration of points with multiplicities in the unit sphere
$S^2\subset\R^3$ can be moved by an element of $SL(2,\C)$ to have centre
of mass the origin if and only if either each multiplicity is strictly less
than half the total, or there are only 2 points and both have the same multiplicity.
\end{thm}

The first case is the stable case, the second the polystable case with a
$\C^*$-stabiliser.

\sec{Example: Grassmannians from GIT and symplectic reduction.} We have seen
how to get $\PP^n$ by GIT and symplectic reduction; we can do something similar
for Grassmannians.

Consider $SL(r,\C)$ acting on $\Hom(\C^r,\C^n),\,r<n,$ linearising
the induced action on the projectivisation $\PP$ of this vector space (we
choose the left action of multiplying on the right by $g^{-1}$).

\begin{prop} \label{grass} $[A]\in\PP$ is stable if $A\in\Hom(\C^r,\C^n)$ has full rank $r$, and unstable otherwise.
\end{prop}

\begin{proof} If rank$(A)<r$ then we can pick a splitting $\C^r=\langle v\rangle\oplus W$ with $A(v)=0$. Then the 1-PS that acts as $\lambda^{r-1}$ on $v$
and $\lambda^{-1}$ on $W$ fixes $[A]\in\PP$ and acts on the line $\C.A$ with
weight $+1$. Therefore $[A]$ is unstable by the Hilbert-Mumford criterion.

Conversely, if $A$ has full rank then, up to the action of $SL(r,\C)$
some multiple of it is the inclusion of the first factor of some splitting $\C^n\cong\C^r\oplus\C^{n-r}$. Diagonalising a given 1-PS, we may assume
further that in this basis we have the action
$$
\text{diag}(\lambda^{\rho_1},\ldots,\lambda^{\rho_r}), \qquad \rho_1\ge\rho_2
\ge\ldots\ge\rho_r, \quad \sum_i\rho_i=0.
$$
Ignoring the trivial 1-PS, there is some $p$ such that $\rho_1=\rho_p>\rho_{p+1}$.
Then the limit $[A_0]$ of $[A]$ under this 1-PS is the inclusion of $\C^p$ as the first $p$ basis vectors of $\C^n$, with the 1-PS acting with
weight $-\rho_1<0$ on $\C.A_0$. Therefore $A$ is stable.
\end{proof}

So the points of the GIT quotient are the injections of $\C^r$ into $\C^n$
modulo the automorphisms of $\C^r$; i.e. they are the \emph{images} of the injections -- the Grassmannian $Gr(r,n)$ of $r$ dimensional subspaces of $\C^n$.

For symplectic reduction, it is easier to consider the affine case of $U(r)<
GL(r,\C)$ acting on $\Hom(\C^r,\C^n)$, with all vector spaces endowed with
their standard metrics. (Above, by working with $\PP$, we had already
divided out by the centre of $GL(r,\C)$ but didn't describe it this way because, as we have seen, it is easier to deal with the linearisation issues in the symplectic picture, where it just amounts to changing the moment map by a central scalar.) The moment map is
\beq{grmm}
A\mapsto i(A^*A-\id),
\eeq
with zeros the orthogonal linear maps that embed $\C^r$ isometrically. Thus
Kempf-Ness recovers the obvious fact that a linear map is congruent by $GL(r,\C)$ to an isometric embedding if and only if it is injective. Dividing
these isometric embeddings by $U(r)$ gives $Gr(r,n)$ again.

\sec{More affine examples.}
Our simple example (\ref{hyp}) has moment map
$$
(|x|^2-|y|^2)/2,
$$
whose zero set intersects each good orbit $xy=\alpha\ne0$ in a unique $U(1)$ orbit $\sqrt\alpha(e^{i\theta},e^{-i\theta})$. It intersects the origin (another
$U(1)$ orbit, corresponding to $\alpha=0$) and misses the other two orbits (the punctured $x$- and $y$-axes). Therefore the symplectic quotient is a copy of $\C$ parameterised by $\alpha$, representing the closed, \emph{polystable} orbits, as anticipated.

If we chose the moment map $(|x|^2-|y|^2+a)/2,\ a>0,$ then we miss the $x$-axis
and the origin, and gain a unique $U(1)$ orbit on the $y$-axis. So the symplectic
quotient is isomorphic, but with a different interpretation. This corresponds
in GIT to a different linearisation, in which the $x$-axis and the origin
are unstable and the punctured $y$-axis is stable. (So this nonclosed orbit
becomes closed upstairs in the new linearisation, and is closed in the locus
of semistable points.)
\medskip

Another standard example is to consider $n\times n$ complex matrices acted
on by the adjoint action of $SL(n,\C)$. The invariant polynomials are the
symmetric functions in the eigenvalues of the matrix (by the denseness of
the set of diagonalisable matrices) -- i.e. functions in the coefficients of the characteristic
polynomial. This reflects the fact that the matrices with
nondiagonal Jordan canonical form have the corresponding diagonal matrices in the
\emph{closure} of their orbits -- all matrices are semistable for this
linearisation (the constant 1 does not vanish on any orbit!),
with the diagonalisable matrices being polystable (their stabiliser is at
least $(\C^*)^n$, after all).

The moment map (for the standard symplectic structure inherited from $\C^{n^2}$)
for the induced action of $SU(n)$ is $A\mapsto\frac12[A,A^*]$ with zeros
the \emph{normal} matrices. Since normal matrices are those that can be orthogonally
diagonalised, the symplectic quotient \{normal matrices\}$\big/SU(n)$ is
the set of diagonal matrices up to the action of the symmetric group, and
so equal to the GIT quotient. (So in this case Kempf-Ness is the obvious fact that a matrix can be diagonalised if and only if it is similar to a matrix that can be orthogonally diagonalised.)

\sec{Back to the Hilbert-Mumford criterion.}
For simplicity of exposition we used the Hilbert-Mumford criterion to prove the Kempf-Ness theorem, to reduce everything to single hamiltonians.
But as we noted there, we could have avoided this and proved it directly by noting that $\log||g.\tilde x||$ is convex on the whole of $G/K$, so an orbit is closed if and only if this log-norm functional is proper, in which case it has a minimum, at which point $m=0$ by (\ref{dd}).

We can then use this to go back and give a sketch proof (more of a discussion, really) of the Hilbert-Mumford criterion. That is we want to show that properness
is equivalent to properness on 1-PSs. As usual one direction
is trivial; for the other one can try to work on $G/K$ as in, for instance, \cite{DK}. The idea is that while 1-PSs cover very little of $G$, since $K$ preserves the norm functional it descends to $G/K$, in which 1-PSs are dense (see the torus case below
where the 1-PSs correspond to directions in $\g/\k\cong\k$ of rational slope). Although it is not a priori clear that properness down each such rational direction is enough to give properness on all of $G/K$, it is clear by openness that if a $G$-orbit is strictly unstable then there will be a rational direction (1-PS) that detects it. So we see that (semi)stability of each 1-PS implies
semistability for $G$.

So this leaves the hard part -- that strict stability for each 1-PS implies
strict stability for $G$. That is, we want to show that if a $G$-orbit is strictly semistable, then there is a 1-PS with zero weight; i.e. the non-properness is detected by a rational direction.

We first show this for $G$ a torus $T^c=(\C^*)^r$.
A $T^c$-action on a vector space splits it into a sum of weight spaces $W_m,\ m\in\mathfrak t^*$, on which $\exp(v)\in T^c,\ v\in\mathfrak t^c,$ acts as the character $\exp(i\langle m,v\rangle)$.

Given any vector $\tilde x$, we let $\Delta_{\tilde
x}\subset\mathfrak t^*$ denote the convex hull of only those weights $m$ in whose weight spaces $\tilde x$ has nonzero components (i.e. the projection
of $\tilde x$ to $W_m$ is nonzero). Any 1-PS corresponds to an integral vector $v\in\mathfrak t$ and so a hyperplane $H_v\le\mathfrak t^*$. The points of $\Delta_{\tilde
x}$ on the negative side of this hyperplane correspond to negative weights in whose weight space $\tilde x$ has a nonzero component, so their existence implies that $\lambda.\tilde x\to\infty$ as $\lambda\to0$ under this 1-PS,
as in (\ref{-ve}). Similarly the existence of points in $\Delta_{\tilde
x}$ on the positive side of the hyperplane prove that $\lambda.\tilde x\to\infty$ as $\lambda\to\infty$.

Thus $\C^*.\tilde x$ is closed, and $\tilde x$ is stable for this 1-PS, if and only if its hyperplane
$H_v\le\mathfrak t^*$ cuts $\Delta_{\tilde x}$ through its interior. Applying this to all integral points $v\in\mathfrak t$ (including those whose hyperplanes
$H_v$ are parallel to the faces of $\Delta_{\tilde x}$) gives the first part
of the following result, which was
explained to me by G\'abor Sz\'ekelyhidi \cite{Sz2}.

\begin{thm} \label{Gabor}
The point $\tilde x$ is stable for every 1-PS if and only if $0\in\mathfrak t^*$ is in the interior of $\Delta_{\tilde x}$, if and only if $\tilde x$ is stable for $T^c$.
\end{thm}

\smallskip \begin{center} \input{interior.pstex_t} \end{center} \smallskip

For the second result we cover the whole of $T^c/T$ (where $T\le T^c$ is the maximal
compact subgroup $T\cong(S^1)^r$) by going in nonrational directions
$v\in\mathfrak t$ too. And if the origin is in the interior, any such $v$ has negative pairing with
at least one of the weights, so the associated orbit (of an analytic
$\C$-subgroup of $T^c$, if $v$ is irrational) will go to infinity as we move
along $v$. In fact the log-norm
function will be proper with a growth that can be bounded below by the minimal
$|\langle m,v\rangle|$. Thus the functional will be proper on all of $T^c/T$,
and we can deduce the Hilbert-Mumford criterion.

Perhaps an easier proof, using the full Kempf-Ness theorem, comes from observing
that the interior of $\Delta_{\tilde x}$ is the moment polytope of the orbit -- the image of $T^c\tilde x$ under the moment map -- so the moment map has a zero in this orbit if and only if the
origin is in the interior.

But the first proof illustrates the key point, which is that the faces of $\Delta_{\tilde x}$ are \emph{rational} -- parallel to hyperplanes $H_v$. So if there is an irrational $v\in\mathfrak t$ that destabilises (has weight $\ge0$)
then since it cannot be contained in a face there is in fact a rational $v$, giving rise to a 1-PS, with the same property. 
So we can avoid the situation of a sequence of stable 1-PSs of negative weight converging to an irrational ``semistable" direction of weight zero lying in a face, making the $T^c$-orbit non-proper but without a 1-PS or rational direction to detect it. \medskip

For an arbitrary group $G$, we can try to reduce to the torus case by dividing $G$ by $K$ on both the left \emph{and} right instead of considering just $G/K$. That is, by spectral theory we can write $G=KT^cK$ for any maximal torus $T^c\le G$; then since the norm functional is invariant under the left hand action of $K$ we are left with proving its properness on
a \emph{compact} family of $T^c$-actions -- the conjugates of the original
action by all $k\in K$.
The result is then basically routine, the point being that in a compact family
of polytopes each containing the origin in its interior, the distance of
the origin to the boundary is bounded below by some $\epsilon>0$.
\medskip

As an application of Theorem \ref{Gabor}, we can strengthen (\ref{smooth})
to recover standard results \cite{GIT, Mu} about which hypersurfaces of degree
$d$ in $\PP^n$ are stable. Namely, forming the Newton polygon of degree $d$ homogeneous polynomials in $(n+1)$ variables, a hypersurface $(f=0)$
defines a subset of integral points of this polytope -- those that appear
in $f$ with nonzero coefficient. Then $(f=0)$ is semistable (or stable) if
and only if these points do not lie to one side of (or strictly to one side
of) any hyperplane
through the centre of the Newton polytope.

\section{Moduli of polarised algebraic varieties $(X,L)$}

\noindent\textbf{The GIT problem.}
This section is unnecessarily technical, and the squeamish reader can skip
it once it is clear why forming moduli of algebraic varieties should
be a GIT problem.

Suppose we want to form a moduli space of polarised algebraic varieties \cite{Mu}.
The polarisation allows us to embed $X$ into a projective space
$$
X\into\PP(H^0(X,L^r)^*),\quad r\gg0.
$$
In fact for $X$ smooth, a theorem of Matsusaka tells us that $r$ can be chosen uniformly amongst all $(X,L)$ with the same Hilbert polynomial $\P(r)=\chi(X,L^r)$. Moreover we can also assume that all higher cohomology groups $H^{\ge1}(X,L^r)$ vanish so that $H^0(X,L^r)$ has dimension $\P(r)$, and that any two $(X_i,L_i)$ are isomorphic if and only if their
embeddings $X_i\into\PP^N,\ N+1=\P(r)$, differ by a projective linear map.

Picking an isomorphism
\beq{iso}
H^0(X,L^r)\cong\C^{N+1},
\eeq
$(X,L)$ defines a point in the \emph{Hilbert scheme} of subvarieties (in
fact subschemes) of $\PP^N$. This moduli space is easy to construct; for
instance as a subscheme of a Grassmannian of subspaces of $S^k(\C^{N+1})^*$;
$X\subset\PP^N$ corresponding to the subspace $H^0(\PP^N,\I_X(k))<
H^0(\PP^N,\O(k))=S^k(\C^{N+1})^*$ of degree $k$ polynomials vanishing on
$X$. The natural Pl\"ucker line bundle then pulls back to give an
anti-ample line bundle on Hilb whose fibre at a point $(X,L)$ is
\beq{line}
\Lambda^{\mathrm{max}}H^0(X,L^{rk})^*\otimes\Lambda^{\mathrm{max}}
S^kH^0(X,L^r).
\eeq
Then we must divide out the choice of isomorphism (\ref{iso}),
i.e. take the GIT quotient of Hilb by $SL(N+1,\C)$.

So by abstract GIT, any choice of $SL(N+1,\C)$-equivariant (anti-)ample line bundle on Hilb gives rise to a notion of \textbf{stability for $(X,L)$}. There are many such; we describe some of those whose associated weights can all be characterised in terms of weights on the line (\ref{line}).

The Hilbert-Mumford criterion requires us to consider $\C^*< SL(N+1,\C)$ orbits of $X\subset\PP^N$. This gives rise to a $\C^*$-equivariant flat family, or \textbf{test configuration}, $(\curly X,\L)\to\C$.

\bigskip \begin{center}
\input{smallfamily.pstex_t}
\end{center} \bigskip

The weight $w_{r,k}$ of the $\C^*$-action on (\ref{line}) is
\beq{wrk}
w_{r,k}=a_{n+1}(r)k^{n+1}+a_n(r)k^n+\ldots,
\eeq
where
$$
a_i(r)=a_{in}r^n+a_{i,n-1}r^{n-1}+\ldots.
$$
Then doing GIT on Hilb with the line (\ref{line}), Mumford's Chow line, or Tian's CM line, gives rise to Hilbert-Mumford criteria that
$\C^*< SL(N+1,\C)$ destabilises $(X,L)$ if $w_{r,k}\succ0$
in the following senses:
\emph{\begin{itemize}
\item HM($r$)-unstable: $w_{r,k}>0$ for all $k\gg0$,
\item Asymptotically HM-unstable: for all $r\gg0,\ w_{r,k}>0$ for all $k\gg0$,
\item Chow($r$)-unstable: leading $k^{n+1}$-coefficient $a_{n+1}(r)>0$,
\item Asymptotically Chow unstable:  $a_{n+1}(r)>0$ for $r\gg0$,
\item K-unstable: leading coefficient $a_{n+1,n}>0$.
\end{itemize}}

To make ``if" into ``iff" requires a few technicalities on the size of $r$;
see \cite{RT1}. In particular K-stability, which
is Donaldson's refinement of Tian's original notion, requires one to pick
a test configuration first, and then choose $r\gg0$. The coefficient $a_{n+1,n}$ is Donaldson's
version of the Futaki invariant of the $\C^*$-action on $(\X_0,L)$; see (\ref{fut}).

There are also notions semistability and polystability in all of these cases;
both defined by nonstrict inequalities, the latter requiring also that
whenever the inequality is an equality, the test configuration should arise
from an automorphism of $(X,L)$, i.e. it should be isomorphic as a scheme
to the product $X\times\C$, but with a nontrivial $\C^*$-action.

In particular we have the following implications (see \cite{RT1}, where our
$a_i$ are denoted $-e_i$): \medskip

Asymptotically Chow stable $\Rightarrow$ Asymptotically Hilbert stable
$\Rightarrow$ Asymptotically Hilbert semistable $\Rightarrow$ Asymptotically
Chow semistable $\Rightarrow$ K-semistable. \smallskip

The increasing number of test configurations that have to be tested as $r\to\infty$
currently prevents one from proving that K-stability implies asymptotic Chow
stability.

\sec{The moment map problem.}
Fix a metric on $\C^{N+1}$ and so $g\_{FS}$ on $\PP^N$ and an induced hermitian
metric on $\O(-1)$. This induces the symplectic
form $\omega_{FS}$ on a smooth $X\subset\PP^N$. This induces a natural symplectic,
in fact K\"ahler, structure on (any smooth subset of smooth points of) Hilb:
\beq{Omega}
\Omega(v_1,v_2):=\int_X\omega_{FS}(v_1,v_2)\frac{\omega_{FS}^n}{n!}\,,
\eeq
where the $v_i$ are the normal components of holomorphic vector fields along $X\subset\PP^N$. This is also (a multiple of) the first Chern class of a natural line bundle on Hilb coming from the ``Deligne pairing" of $\O_X(1)$
with itself $(n+1)$-times \cite{Zh}.

Let $m\colon\PP^N\into\su(N+1)^*$ denote the usual moment map (\ref{mmpn}).

Then \cite{Do3}, just as for a finite number of points in $\PP^1$ (\ref{mm1}), the moment map for $SU(N+1)\acts(\Hilb,\Omega)$ takes
$X\subset\PP^N$ to a multiple of its centre of mass in $\su(N+1)^*$:
\beq{mmb}
\mu(X)=\int_X m\,\frac{\omega_{FS}^n}{n!}\ \in\ \su(N+1)^*.
\eeq

So zeros of moment map correspond to \textbf{balanced} varieties $X\subset\PP^N$.

The fact that Hilb is not smooth means there are complications in applying
the
Kempf-Ness theorem directly, but nonetheless the following is an essentially
finite dimensional result. It was first proved by Zhang \cite{Zh}, and then rediscovered and reproved in different forms by Luo, Paul, Wang and Phong-Sturm.

\begin{thm} \label{KN}
$X\subset\PP^N$ can be balanced by an element of $SL(N+1,\C)\ \Longleftrightarrow\ X$ is Chow polystable.
\end{thm}

The balanced condition can be re-expressed as follows. The metric on $\O_X(1)=L^r$ is the quotient
metric induced from that on $H^0(\O_X(1))=H^0(X,L^r)$ by the surjection of
vector bundles
$$
\underline{H^0(\O_X(1))}\to\O_X(1)\to0
$$
on $X$. So picking an orthonormal basis $\sigma_i\in H^0(X,L^r)$, we have the identity
\beq{e1}
\sum_i|\sigma_i(x)|^2\equiv1
\eeq
on $X$. (More generally, given an orthonormal basis $e_i$ of an inner product space $V$ and a surjection $V\rt{\pi}W$, we have
the identity $\sum_i|\pi(e_i)|^2=\dim W$ in the induced metric on $W$. Given
any basis $\sigma_i$ of $H^0(L^r)$ the above expression (\ref{e1}) is the
pointwise ratio of the given metric
on $L^r$ and the Fubini-Study metric on $L^r$ induced by embedding in
$H^0(L^r)^*$ and pulling back the metric gotten by declaring the $\sigma_i$ to be orthonormal. This is constant if and only if the metric on $(X,L^r)$ really is such a Fubini-Study metric.)

But then in these coordinates, the moment map (\ref{mmb}) constructed using
(\ref{mmpn}), takes $X$ to the matrix with $(ij)$th entry
\beq{bal2}
i\left(\int_X\sigma_i(x)\sigma_j(x)^*\frac{\omega_{FS}^n}{n!}-
\frac1{N+1}\delta_{ij}\right)\in\su(N+1).
\eeq
Thus the balanced condition is equivalent to the $\sigma_i$ being orthonormal
(up to a constant scale) in the induced $L^2$-metric on $H^0(\O_X(1))$.
That is, up to scale, the original metric on $\C^{N+1}\cong H^0(\O_X(1))=
H^0(X,L^r)$ equals the $L^2$-metric given by integration against $g\_{FS}|\_X$. By (\ref{e1}) this is equivalent to 
\beq{bal3}
\sum_i|s_i(x)|^2\equiv\mathrm{const},
\eeq
where the $s_i$ are now an orthonormal basis \emph{with respect to the $L^2$-metric on $H^0(X,L^r)$} (rather than the original metric).
A final way of saying this is that starting with a metric on $\C^{N+1}$ we can induce another by first inducing the Fubini-Study metric on $X\subset\PP^N$
and the hermitian metric on $\O_X(-1)$, and then using this to give, by integration,
an $L^2$-metric on $\C^{N+1}=H^0(\O_X(1))$. Balanced metrics are then the
fixed points of this operator.

\sec{Asymptotics as $r\to\infty$.}
Fix a metric on $(X,L)$ (e.g. by picking a metric on $H^0(L)$ and then
inducing the Fubini-Study metric on $X\subset\PP(H^0(L^*))$ and $L=\O(1)|_X$).
This then induces one on $L^r$ for all $r$, and so $L^2$-metrics on $H^0(X,L^r)$
for all $r$.

Picking an $L^2$-orthonormal basis $s_i\in H^0(X,L^r)$, we can then define, for each $r$, the \textbf{Bergman kernel}
\beq{Ber}
B_r(x_1,x_2)=\sum_is_i(x_1)\otimes s_i(x_2)^*
\eeq
on $X\times X$. This is the integral kernel for the $L^2$-orthogonal projection
of $C^\infty$ sections of $L^r$ onto holomorphic sections. Restricting
to the diagonal gives
\beq{Berg}
B_r(x):=\sum_i|s_i(x)|^2,
\eeq
so the balanced condition (\ref{bal3}) is equivalent to $B_r$ (\ref{Berg})
being \emph{constant on $X$}.

This expresses the finite dimensional balanced condition (a condition
for a metric on $H^0(X,L^r)$) as a \emph{pointwise} condition for a metric
on $(X,L)$ (a fact that will be explained later via Donaldson's
double quotient construction) and we can look at the asymptotics of the ``density
of states" function $B_r(x)$ as $r\to\infty$ and expect it to only depend
on local differential-geometric data. This is because, as is well
known to quantum physicists and is made precise in \cite{Ti3}, one can form a basis of sections of $H^0(L^r)$ whose norms are approximately peaked Gaussians concentrated in balls of radius const$/\sqrt{r}$ and so volume const$/r^n$.
(These are the \emph{coherent states} of geometric quantisation; under the
metric isomorphism $H^0(L^r)\cong H^0(L^r)^*$ they correspond to evaluating
sections at points -- the centres of the peaks.)
The relationship between the volume of small balls about $x\in X$ and the scalar curvature $s(x)$ at $x$ means
that as $r\to\infty$ the number of peaked sections that can be packed into a fixed ball of volume $\epsilon$ about $x$ is $\sim\epsilon
(r^n+\frac n2s(x)r^{n-1}+\ldots)$. Globally this gives rise to
$$
\vol(X)r^n+\frac n2\int_X\!s\,\frac{\omega^n}{n!}r^{n-1}+O(r^{n-2})\ =\ 
\frac{L^n}{n!}r^n-\frac{K_X.L^{n-1}}{2(n-1)!}r^{n-1}+O(r^{n-2})
$$
sections -- approximating the Riemann-Roch formula.

In fact, as $r\to\infty$ ($\Rightarrow N\to\infty$) $B_r(x)$ has an asymptotic expansion (Tian, Zelditch, Catlin, W.-D. Ruan, Z. Lu)
\beq{expansion}
B_r(x)\sim r^n+\frac1{2\pi}s(x)r^{n-1}+O(r^{n-2}).
\eeq
More precisely, $||B_r(x)-r^n+\frac1{2\pi}s(x)r^{n-1}||\_{C^\alpha}\le
Cr^{n-2}$ for $\alpha\ge0$, where the constant $C$ depends on both $\alpha$
\emph{and} the metric -- it can only be taken to be uniform for metrics in
a \emph{compact subset} of the space of metrics.

Roughly speaking then, balanced metrics should tend towards cscK metrics with $[\omega]=[c_1(L)]$. What we have seen so far should motivate the following
results.

\begin{thm} \label{thm}
If $(X,L)$ admits a cscK metric in $[c_1(L)]$ and has finite automorphism group then $(X,L^r)$ can be balanced in $\PP(H^0(X,L^r)^*)$ for $r\gg0$. Thus it is Chow stable, and so K-semistable. \\
The metrics given
by $r^{-1}$ times by the pull backs of the balanced metrics converge to the cscK metric. \smallskip \\
Conversely if $(X,L^r)\subset\PP^{N(r)}$ is balanced for $r\gg0$ \emph{and} the resulting $\omega\_{FS,r}$ are convergent, then the limit metric is cscK.
\smallskip \\
Finally, cscK metrics compatible with a fixed complex structure are unique up to holomorphic automorphisms of $X$.
\end{thm}

This result is due to Donaldson \cite{Do3}; we will discuss the proof
of the balanced result in a later section. Tian had previously proved K-semistability
for KE metrics \cite{Ti2}, and a related convergence result for sequences
of Fubini-Study metrics \cite{Ti3}, following a suggestion of Yau \cite{Y1}. Using \cite{Do3} Mabuchi proved that cscK manifolds with automorphisms
are Chow polystable if the automorphisms satisfy a certain stability condition
\cite{Mb2}. Donaldson \cite{Do5} then showed that cscK $\Rightarrow$ K-semistable
without any condition on automorphisms.

Uniqueness was originally proved by Bando-Mabuchi \cite{BM} for KE metrics,
by Chen \cite{Ch} for cscK metrics when $c_1\le0$, then by Donaldson in the general cscK case with finite automorphisms. Again the finite automorphisms condition was relaxed by Mabuchi, and, in the more general setting of extremal metrics and K\"ahler non-projective metrics, Chen-Tian \cite{CT}.

When our polarisation $L$ is a power of the canonical bundle $K_X$, then cscK metrics are in fact KE: those with Ricci form (the induced curvature of $K^{-1}_X$) a constant multiple of the \K form. It is clear that KE metrics are cscK; the converse follows from the calculation
$$
\Delta\,\text{Ric}=-i\partial\overline\partial\Lambda\,\text{Ric}=
-i\partial\overline\partial s.
$$
Here $\Lambda$ is the adjoint of $\omega\wedge$, and we use the fact that
Ric is a closed real $(1,1)$-form, so it is $\partial$- and
$\overline\partial$-closed. $s=\Lambda\,$Ric is the
scalar curvature. So for $s$ a constant, Ric is harmonic, as is $\omega$. But, after scaling, they represent the same cohomology class, and so are identically
equal. \medskip

KE metrics were first proved to exist on compact \K manifolds with positive canonical bundle by Aubin \cite{Au} and Yau \cite{Y}, and with trivial canonical bundle by Yau \cite{Y}.
It was Calabi \cite{Ca} who initiated the study of cscK and \emph{extremal} metrics:
those which extremise the Calabi functional $\int_Xs^2\omega^n$ over cohomologous \K forms; they are the metrics with $\nabla s$ a holomorphic vector field. 
Apart from Aubin and Yau's (nonconstructive) results, there are few compact
examples of cscK or KE metrics. Siu \cite{S}, Tian and Nadel \cite{Na} found examples with symmetry,
Tian showed Fano surfaces with reductive automorphism
groups admit KE metrics \cite{Ti1}, Burns-de Bartolomeis \cite{BdB} and Hong
\cite{Ho} gave constructions of cscK metrics on certain projective bundles
over cscK bases, and there are constructions for blow ups of these \cite{AP,
LB, RS} and smooth fibrations of cscK manifolds \cite{Fi}. Bourguignon \cite{Bo}
and Biquard \cite{Bi} have given excellent surveys of KE and cscK metrics respectively.

\sec{An example -- blow ups of cscK manifolds.}
The results of \cite{AP} give a beautiful illustration of the theory described here and the link between cscK and balanced metrics. Arezzo and Pacard consider a cscK manifold $(X,\omega)$ and its blow up in some points $p_i$,
\begin{eqnarray*}
\pi\colon\widehat X & \To & X, \\
E_i & \to & p_i.
\end{eqnarray*}
It is proved that there is a cscK metric in the class $\pi^*[\omega]-\epsilon\sum_im_i[E_i]$
for $\epsilon,m_i>0$ and $\epsilon$ sufficiently small, if the $m_ip_i$ satisfy two conditions with respect to Aut$(X,\omega)$. Arezzo and Michael Singer
observed that one of these conditions could be rewritten as a balanced condition.
Namely there is a moment map $X\Rt{\mu_a}\mathfrak{ham}(X,J,\omega)^*$ for the action of the hamiltonian isometry group of $X$, and the conditions are that
\beq{MAP}
\sum_im_i\mu_a(p_i)=0, \text{ and the } \mu_a(p_i) \text{ generate }
\mathfrak{ham}(X,J,\omega)^*.
\eeq

We can interpret this in the projective case, where $\mathfrak{ham}(X,J,\omega)$
becomes $\mathfrak{aut}(X,L)$, as follows. Taking $\epsilon$ very small is equivalent to replacing
the polarisation by a very large power $r\gg0$, whereupon the cscK condition approximates the balanced condition (\ref{expansion}) (for what follows we only need that the approximation is valid for the linearisation of the equations
as $r\to\infty$). Then morally,
in replacing $(X,L^r)$ by $(\widehat X,\pi^*L^r(-\sum m_iE_i))$ we are perturbing
a balanced $X\subset\PP^N=\PP(H^0(L^r)^*)$ only a little bit and so end up
with a manifold that is nearly balanced. Slightly more precisely, set $I=
H^0(L^r\otimes\I_{\cup_im_ip_i})$ and split the exact sequence
\spreaddiagramcolumns{1pc}
$$
0\To H^0(L^r\otimes\I_{\cup_im_ip_i})\To \xymatrix{H^0(L^r) \rto
\ar@/^.9pc/@{<--}[r] & \bigoplus_i\C_{p_i}^{m_i}} \To0
$$
by picking peaked approximately Gaussian sections of $L^r$ on $X$ at the $p_i$, as in our discussion of (\ref{expansion}).
Away from the $p_i$, therefore, points in the image of $X\into\PP(H^0(L^r)^*)$
almost annihilate this $\bigoplus_i\C_{p_i}^{m_i}$, i.e. they lie very close to
$\PP(I^*)$, as in the following diagram.

\begin{center} \input{balanced.pstex_t} \end{center}
\spreaddiagramcolumns{-1pc} \smallskip

The dashed arrows denote the rational map $\xymatrix{\PP(H^0(L^r)^*) \ar@{-->}[r] & \PP(I^*)}$ blowing up $\PP(\bigoplus_i\C_{p_i}^*)$;
on restriction to $X$ this blows up the $p_i$ and embeds the result in $\PP(I^*)$.

So away from the
$p_i$, the moment map $\PP^N\to\su(N+1)^*$ of (\ref{mmpn}), projected to $\su(I)^*$, is very close (as $r\to\infty$) to the rational projection to
$\PP(I^*)$ followed by the moment map $\PP(I^*)\to\su(I)^*$.

Since the exceptional divisors are small, we can integrate over $X$ (or its
blow up in the $p_i$) to find that the centre of mass in $\su(I)^*$ is close
to the projection of that in $\su(N+1)^*$. But this is zero, so $\widehat
X$ is close to balanced, as claimed.

Now the exact sequence expressing the derivative $D$ of the $SU(N+1)$ action
on the Hilbert scheme of $\PP(H^0(X,L^r)^*)\supset X$,
$$
0\to\mathfrak{aut}(X,L^r)\to\su(N+1)\Rt{D}T\Hilb\cong T^*\Hilb
$$
(with the last isomorphism induced by the symplectic form), has dual
\beq{abc}
0\ot\mathfrak{aut}(X,L^r)^*\ot\su(N+1)^*\stackrel{d\mu}{\longleftarrow}T\Hilb,
\eeq
by the definition of the moment map $\mu=\int_Xm\,\omega_{FS}/n!$
(\ref{mmb}).

If the automorphism group of $(X,L^r)$ is finite (so the condition (\ref{MAP}) is vacuous) then $D$ is injective and
its adjoint $d\mu$ is \emph{onto}. So we expect to be able to move a little in the orbit to move back to a balanced metric with $\mu=0$ to correct the
perturbation introduced by the $p_i$. This of course involves some estimates, which
is what \cite{AP} work out for the cscK problem, to show that for $\mathfrak{aut}=0$
there is always a cscK metric on the blow up.

When the automorphism group is nontrivial this map $d\mu$ is \emph{not} onto, so we must ensure that on perturbing as above we end up inside its image to apply the same
argument. That is, by (\ref{abc}), the image of the moment map in
$\mathfrak{aut}(X,L^r)^*$ should be zero. Since the moment map is the centre of mass, and since
we have added masses $m_i$ at the exceptional divisors $E_i$ lying over $p_i$, we must ensure that, to first order, the $\cup_im_ip_i$ should be balanced
in $\mathfrak{aut}(X,L^r)^*$.

This recovers (\ref{MAP}) as the necessary linearised condition. The second condition
is a nondegeneracy condition that allows one to perturb the metric on and
around the exceptional divisors to move the moment map enough to solve the
equation to higher orders. 

As pointed out by Donaldson, Hong's results \cite{Ho} on when a cscK metric exists on the projectivisation
of a HYM bundle over a cscK base involves a similar moment map condition
for the action of the automorphism group of the base on the moduli of vector bundles. \medskip

These examples illustrate a general principle about moment map problems: that transverse
(regular) points of $\mu^{-1}(0)$ have no automorphisms, whereas for nontransverse
points $x$ the cokernel of $d\mu$ is canonically $(\mathfrak g^x)^*$, the dual of the Lie algebra of the stabiliser subgroup of the point $x\in X$. Thus when one perturbs a solution $x$ of $\mu=0$ with stabiliser subgroup $G^x<G$, the obstruction to extending a first order deformation
lies in $(\mathfrak g^x)^*$, and is nothing but the derivative of the moment map of the action of $G^x<G$. This follows from the exact sequence
$$
T_xX\Rt{d\mu}\mathfrak g^*\to(\mathfrak g^x)^*\to0,
$$
the dual of
$$
0\to\mathfrak g^x\to\mathfrak g\to T^*_xX,
$$
with the last map the composition of the $\mathfrak g$-action on $T_xX$ and
contraction with the symplectic form (cf. (\ref{abc})).

\sec{The infinite dimensional setup.}
Instead of letting the dimension $N$ of our quotient problem go to infinity, Donaldson
\cite{Do1} also gave a purely infinite dimensional formal symplectic quotient formulation.

The group of Hamiltonian diffeomorphisms acts on $(X,\omega)$ and so on the
space of complex structures which make $(X,\omega)$ K\"ahler:
$$
\text{Ham}(X,\omega)\acts\mathcal J:=\{\text{$\omega$-compatible complex structures on } X\}.
$$
Acting by pullback, the infinitesimal action of a hamiltonian $h$, with hamiltonian
vector field $X_h$, on a complex structure $J$ is $\mathcal L_{X_h}J$. At
the Lie algebra level this can be complexified so that $ih$ acts as
$$
J\mathcal L_{X_h}J=\mathcal L_{JX_h}J=\mathcal L_{X_{ih}}J,
$$
by the integrability of $J$. Thus it acts through the vector field
$$
X_{ih}:=JX_h.
$$
We note that the action of this vector field on $\omega$ is
$$
\mathcal L_{JX_h}\omega=d(JX_h\ip\omega)=d(Jdh)=d(-i\partial h+i\dbar h)=
2i\partial\dbar h,
$$
changing $\omega$ within its cohomology class by the \K potential $h$
to another form compatible with $J$.

We can contract these vector fields with $\omega$ to write them as one-forms. By Hodge theory,
$$
\Omega^1(X)=dC^\infty(X)\ \oplus\ H^1(X,\R)\ \oplus\ d^*\Omega^2.
$$
The first summand corresponds to the hamiltonian vector fields, the second
to symplectomorphisms modulo those which are hamiltonian, and inside the
third lies $d^*(C^\infty(X)\omega)$ as those which preserve the compatibility of $\omega$ with $J$ (i.e. down which the Lie derivative of $\omega$ is of
type (1,1)).
These constitute the complexified hamiltonian action, by the \K identity
$$
d^*(h\omega)=i(\dbar-\partial)h=Jdh,
$$
whose contraction with (the inverse of) $\omega$ is $JX_h=X_{ih}$.

So, assuming $H^1(X,\R)=0$ for simplicity,
integrating up this complexified Lie algebra suggests defining the complexification of Ham$(X,\omega)$ to be the set of diffeomorphisms of $X$
such that the pullback of $\omega$ is compatible with $J$ (i.e. of type (1,1)):
\beq{Gc}
\big\{f\colon X\to X\ \colon\ \exists\,h\in C^\infty(X,\R) \text{ such that
} f^*\omega=\omega+2i\partial\dbar h\big\}.
\eeq
While this description depends on $J$, it does formally complexify
Ham$(X,\omega)$: we have already seen that it has the right tangent space
$C^\infty(X,\R)_0\otimes\C$ at each point, and it is, crucially, \emph{contractible onto} Ham$(X,\omega)$ by Moser's theorem and the convexity of the space of \K forms.

\sec{The complexified orbits.}
Although (\ref{Gc}) is not actually a group, its \emph{orbits} on
$\mathcal J$ (consisting of pullbacks of complex structures by the above diffeomorphisms) make perfect sense and complexify the Ham$(X,\omega)$ orbits.

Since any two complex structures $J,\,f^*J$ in such an orbit differ by a diffeomorphism, we consider them isomorphic. They are both, by construction,
compatible with $\omega$, but the \K structures $(J,\omega),\,(f^*J,\omega)$
they define need not be isomorphic as the latter is only isomorphic to
$(J,(f^{-1})^*\omega)$. Pulling back by the diffeomorphisms $f$ in this way
(i.e. fixing $J$ and moving $\omega$ instead of the other way round) we get an exact sequence
\begin{multline} \label{corbit}
\text{Ham}(X,\omega)\to\text{Ham}^c(X,\omega).J\twoheadrightarrow \\
\big\{\text{compatible \K metrics on $(X,J)$ in the $H^2$ class $[\omega]$}\big\}.
\end{multline}
The last arrow is \emph{onto} because any such $\omega'$
is of the form $\omega+2i\partial\dbar h$, and so 
diffeomorphic to $\omega$ (since by the convexity of the space of \K forms it is connected to $\omega$ through a family of \K forms $\omega+ti\partial\dbar h$ which are therefore all diffeomorphic by Moser's theorem). Thus the space
of \K metrics on $(X,J)$ is formally of the form $G/K$. This sequence should be compared to its (more familiar) bundle analogue in (\ref{borbit}).

The set-theoretic ``quotient"
by the complexified group (i.e. the set of complexified orbits) is therefore
the set of isomorphism classes of integrable complex structures on $X$ (that
are compatible with one of the symplectic forms $f^*\omega$).

\sec{Moment map = scalar curvature.}
The \K structure on $X$ induces one on $\mathcal J$ by integration.
This is preserved by Ham$(X,\omega)$, and we can ask for a moment map. Considering $C^\infty(X,\omega)_0$ (the functions of integral zero) to lie in the dual of the Lie algebra $C^\infty(X,\R)\big/\R$ by integration against $\omega^n$, and setting $s_0$ to be the topological constant $\int_Xc_1(X).\omega^{n-1}\big/
\int_X\omega^n=\int s\omega^n\big/\int\omega^n$ (the average scalar
curvature), Fujiki \cite{Fj} and Donaldson \cite{Do1} show that
\beq{ms}
\text{Moment map }=s-s_0.
\eeq
This should be no surprise, since we were looking for a function depending algebraically on the second derivatives of the metric, i.e. an invariant scalar derived from the curvature, which can only be a multiple of the scalar
curvature. Thus zeros of the moment map correspond to cscK metrics.

\sec{Norm functional = Mabuchi's K-energy.}
The formula (\ref{imv2}) for the change in the log-norm functional $M=\log
||\tilde x||$ along a complexified orbit, gives the following in this
infinite dimensional set-up. Moving down the orbit of $ih,\,h\in C^\infty(X,\R)$,
i.e. in the family of \K forms $\omega_t=\omega+2it\partial\dbar h$,
\beq{ddt}
\frac{dM}{dt}=m_h,
\eeq
where $m_h=\langle m,h\rangle$ is the hamiltonian function on $\mathcal J$ for the element of the Lie algebra $h\in C^\infty(X,\R)$. Since the moment map $m=s-s_0$
(\ref{ms}), $m_h=\int_X(s-s_0)h\omega_t^n\big/n!$, and
$$
M(\omega_s)=\int_0^s\int_X(s_t-s_0)h\frac{\omega_t^n}{n!}\ dt,
$$
where $s_t$ is the scalar curvature of $\omega_t$. This is precisely the
\emph{Mabuchi functional} or \emph{K-energy} \cite{Mb1}, defined up to a
constant (equivalent to the ambiguity in the choice of a lift of a point to the line bundle above it). It can indeed be written as the log-norm
functional for a Quillen metric on a line bundle over the space of \K metrics;
see for example \cite{MW}.
Its critical points are cscK metrics,
and one expects such a metric to exist on $(X,J)$ if and only if $M$ is proper
on the space of \K metrics on $(X,J)$ (which is the infinite dimensional
analogue of $G/K$ by (\ref{corbit})).

\sec{Weight = Futaki invariant.}
The formula (\ref{ddt}) at a fixed point (e.g. the limit point of a 1-PS
when this exists and is smooth), on the line over which $\C^*$ acts with weight $\rho$, is
\beq{fut}
\frac{dM}{dt}=\rho=m_h=\int_X(s-s_0)h\frac{\omega^n}{n!}\,.
\eeq
Compare (\ref{imv}, \ref{imv2}, \ref{rho}). This is the statement that ``the derivative of the Mabuchi energy is the Futaki invariant" \cite{Mb1, DT}.

The right hand side is, up to a sign, the original definition of the Futaki invariant \cite{Fu}
for a smooth polarised manifold $(X,L)$ with a $\C^*$-action. Noting as above
that it is the weight of the induced action on a line led Donaldson to give the more general definition $a_{n+1,n}$ described earlier, for an arbitrary polarised scheme $(X,L)$.

\sec{Approximation and quantisation.} As Donaldson explains in \cite{Do6}, the finite dimensional problem of balanced metrics can be thought of as the quantisation of the infinite dimensional problem of cscK metrics, which emerges
as the classical limit as $r,N\to\infty$.

As in quantum theory we think of the spaces of sections of the line bundle
$L^r$ as wave functions on $X$, with a basis of Gaussian sections, peaked
around points on $x$. As $r\to\infty$ these peak more, largely supported
in balls of radius const$/r$. Our $SL(N+1,\C)$ group action moves
these sections around the manifold, which may be thought of as moving quantised
chunks of manifold of volume $\sim 1/r^n$ around $X$ (thanks to Anton Gerasimov
for this analogy). In the limit this is meant to approximate the classical
limit of the diffeomorphisms (\ref{Gc}) in the complexification of $\tH$ moving points of the manifold around.

There is in fact a natural map $\su(N+1)\to C^\infty(X,\R)$, though it is
only a homomorphism of Lie algebras to leading order in $r$ \cite{CGR}. A skew-adjoint
endomorphism $iA\in\su(N+1)$ gives an infinitesimal automorphism of $\PP^N$ whose vector field $v_A$ is hamiltonian with respect to the Fubini-Study symplectic form. Its hamiltonian is the function (Berezin symbol)
\beq{symbol}
\PP^N\ni x=[\tilde x]\mapsto\frac{\langle A\tilde x,\tilde x\rangle}
{||\tilde x||^2}=:h_A.
\eeq
On $X$, $h_A|_X$ induces a hamiltonian vector field which is the orthogonal
projection of $v_A$ from $T\PP^N|_X$ to $TX$.

Using the fact that $TX$ is invariant under the complex structure $J$, and
working with complexified hamiltonian vector fields (of the form $X_h+JX_g=:
X_{h+ig}$), the same working shows that the same formula defines a map from $\mathfrak{sl}(N+1,\C)$ to the Lie algebra $C^\infty(X,\C)$ of the complexification (\ref{Gc}) of $\tH$. Thus the change in metric on $X$ induced by
pulling back the metric along an $\mathfrak{sl}(N+1,\C)$ vector field in
$\PP^N$ is the same as that induced by pulling back along its orthogonal
projection \emph{tangent to $X$}.  (Thanks to G\'abor Sz\'ekelyhidi for this observation \cite{Sz2}.)

In this way algebraic 1-PS orbits, i.e. test configurations,
give rise to curves in the complexification of the $\tH$-action
on $\mathcal J$ which approximate 1-PS orbits. Using the description of these orbits in terms of a fixed
complex structure and varying \K form (\ref{corbit}), this simply corresponds to restricting the Fubini-Study metric of $\PP^N$ to the test configuration.

To get a map back we orthogonally project the prequantisation representation
$\tH\to$\,Aut$(\Gamma(L^r))$ to $H^0(L^r)<\Gamma(L^r)$ using the Bergman kernel (\ref{Ber}). That is $h\in C^\infty(X,\R)$ maps to the infinitesimal automorphism $iA\in\su(H^0(L^r))$ defined by
$$
iA(s)=\sum_i\langle\nabla_{X_h}s+ihs,s_i\rangle\_{L^2}s_i.
$$
Again, this is not a homomorphism (except to leading order in $r$). The problem is that we had
to use the Bergman kernel because quantisation is \emph{not} a symplectic
invariant (it cannot be done equivariantly with respect to symplectomorphisms
or elements of $\tH$). That is, it is not independent of choices of complex structure because the pullback of $s\in H^0(L^r)$ by $\tH$ is not in general
holomorphic.

\sec{Donaldson's double quotient construction.}
Because of this problem Donaldson \cite{Do3} considers pairs of a complex
structure $J\in\mathcal J$ and a section $s\in\Gamma(L^r)$ which is \emph{holomorphic
with respect to $J$}; these are clearly acted on by $\tH$. In fact he considers
$N+1=h^0(X,L^r)$-tuples of sections:
$$
\mathcal S=\big\{(J,\{s_i\})\in\mathcal J\times\Gamma(L^r)^{N+1}\ \colon\,
\dbar_Js_i=0,\ i=1,\ldots,N+1\big\}.
$$
Here, as usual, $L$ has a metric and hermitian connection, and $\dbar_J$ is the the $(0,1)$-part of the induced connection on $L^r$ with respect to the complex structure $J$. Since the curvature $2\pi ri\omega$ is compatible with $J$ (by the definition
of $\mathcal J$), i.e. of type (1,1), $\dbar_J^2=0$ and $\dbar_J$ defines an integrable holomorphic structure on $L^r$ by the Newlander-Nirenberg theorem.

We now have actions of $GL(N+1,\C)$ \emph{and} $\tH^c$. These commute, and
both have centre $\C^*$ acting by scalars on $L$, so we can quotient by $\tH^c$
and then $SL(N+1,\C)$, or by $GL(N+1,\C)$ and then Ham($X,\omega)^c$. In
this way we will see how an infinite dimensional moment map problem is equivalent
to a finite dimensional one.

Dividing by $GL(N+1,\C)$ leaves $\mathcal J$ (with a fibration
over it by the Grassmannian of $(N+1)$-planes in $H^0(L^r,\dbar_J)$, by
Proposition \ref{grass}, but for $L$ sufficiently ample $N+1=h^0(L^r)$ and this is a single point).
In turn the formal complex quotient of this by Ham($X,\omega)^c$, discussed above, is the space of complex structures on $X$ (compatible with some symplectic structure in the diffeomorphism group orbit of $\omega$). Taking symplectic
reductions instead we end up with cscK metrics (together with orthonormal
bases of $H^0(L^r)$ modulo the unitary group -- i.e. just a point). So far
then, we have just reproduced what we already knew.

However, we can put a different symplectic structure $\Omega_r$ on $\mathcal J$, and one that tends to $\Omega$ as $r\to\infty$. Namely, the fact that
the $s_i$ determine an embedding of $X$ into $\PP(H^0(L^r)^*)$ for $r\gg0$
means that the natural projection
$$
\mathcal S\to\Gamma(L)^{N+1}
$$
is an embedding, and we can pullback the natural $L^2$-symplectic form from
the latter to define $\Omega_r$.

Now \cite{Do3, Do6} the moment map for the Ham$(X,\omega)$-action becomes
$$
\Big(\frac12\Delta+r\Big)\sum_i|s_i(x)|^2,
$$
with zeros the solutions of $\sum_i|s_i(x)|^2=$\,constant.

If we first take the symplectic reduction by Ham$(X,\omega)$ then this involves solving
$\sum_i|s_i(x)|^2=$\,constant, which we have already observed in (\ref{e1})
says that the
metric on $X$ is the restriction of the Fubini-Study metric on $\PP(H^0(L^r)^*)
\supset X$ when we put the metric on $H^0(L^r)$ that makes the $rs_i$ orthonormal
(the scaling arising because we have ignored the central scalar action).
But since this is a \K metric in the same class as $\omega$, we have already observed (\ref{corbit}) that we can solve this in a Ham$(X,\omega)^c$ orbit, uniquely up to the action of Ham$(X,\omega)$. Next we take the reduction by $SU(N+1,\C)$,
which by (\ref{mmb}) means we try to balance $X\subset\PP(H^0(L^r)^*)$ in
the metric in which the $rs_i$ are orthonormal. By Theorem \ref{KN} there
is a solution to this problem in a $SL(N+1,\C)$-orbit of $X$, unique up to
$SU(N)$, if and only if $X\subset\PP(H^0(L^r)^*)$ is Chow polystable.

So that gives us the finite dimensional problem of solving (\ref{bal2}),
(which, as observed there, is equivalent to the metric on $H^0(L^r)$ being
the $L^2$-metric).

Taking the symplectic reduction in the opposite direction gives instead the pointwise description (\ref{bal3}) of the balanced condition. Namely, first taking the reduction by $SU(N+1)$ gives us an \emph{orthonormal} basis $s_i$ (up to an overall scale which could be removed by putting back the central $\C^*$-action) in each $SL(N+1,\C)$ orbit, unique up to $SU(N+1)$,
if and only the original $s_i$ were
linearly independent (Proposition \ref{grass}). Then taking the reduction
by Ham$(X,\omega)$ involves solving (\ref{bal3}) $B_r(x)=$\,const for the
metric. So we see how solving this infinite dimensional moment map problem
has been reduced to the finite dimensional balanced moment map problem.

This latter equation has the advantage that it is asymptotically close to
the cscK equation (\ref{expansion}). If quantisation really ``worked" it
would be exactly the cscK equation, and proving Donaldson's result that cscK
$\Rightarrow$ balanced would be trivial. Since it is only asymptotically
close, Donaldson crucially uses the ``failure" of quantisation to move from
a cscK solution to a balanced solution, as we now describe.

\sec{CscK $\Rightarrow$ balanced.}
In \cite{Do3}, Donaldson proves a ``quantitative" version of the Kempf-Ness theorem: if the moment map $m(x)$ at a point $x$ is small, and the action of the Lie algebra at $x$ is injective, with a sufficiently large lower bound
on its smallest eigenvalue in a sufficiently large neighbourhood of $x$,
then there exists a zero of $m$ close to $x$ in its complexified orbit.
Flowing down the gradient of $-||m||^2$, i.e. down $JX_{m^*}$ (where $m^*\in\k$ is dual to $m\in\k^*$ under the Killing form), the conditions ensure that
$X_{m^*}$ is sufficiently large and so $||m||^2$ decreases sufficiently fast
for sufficiently long to converge to a zero of $m$.

He applies this to the $SU(N+1)$-action on the symplectic reduction of $\mathcal
S$ by Ham$(X,\omega)$. The cscK metric ensures that we are close to a balanced
metric (zero of the moment map) as $r\to\infty$. Then to give a lower bound
for the injectivity of the $\su(N+1)$-action it is equivalent to give a bound
for the orthogonal projection of its action perpendicular to the orbits of
Ham$(X,\omega)$ upstairs on $\mathcal S$.

Donaldson shows that the projection of the action of $iA\in\su(N+1)$ onto
the tangent to the Ham$(X,\omega)$ orbits is just what one might expect from quantisation: it is the action of its Berezin symbol $h_A$ (\ref{symbol}).
So the normal projection we require is given by the difference in the
actions of $iA$ and $h_A$ on $\mathcal S$.

It is here is where the failure of quantisation to be equivariant with respect
to Ham$(X,\omega)$ is used -- to show that this difference is sufficiently
large in some sense. Of course quantisation is invariant with respect to
\emph{holomorphic} hamiltonian vector fields, i.e. those functions satisfying
$$
\mathcal Dh:=\dbar X_h=\dbar(dh\ip\omega^{-1})=0.
$$
\cite{Do3} assumes that Aut$(X,J)=0$, so that $\ker\mathcal D$ is just
the constants. Then the (fixed) lowest eigenvalue of $\mathcal D$ gives the
lower bound on the difference of the actions of $iA$ and $h_A$. This gives the required estimates, as $r\to\infty$ and we get closer to a zero of the balanced moment map equation, to apply Donaldson's quantitative Kempf-Ness theorem.
\medskip

So $SU(N+1)$ really ``approximates" $\tH$, in the sense that its finite dimensional moment map converges to the infinite dimensional one (\ref{expansion}),
the symplectic structures $\Omega_r\to\Omega$, and the natural norm functionals and weights tend to their infinite dimensional analogues (the Mabuchi functional
and Futaki invariant) as $r\to\infty$; see \cite{Do6}.
Also the space of ``algebraic metrics" (the restrictions of the Fubini-Study
metrics $SL(N+1,\C).\omega_{FS}$ from $\PP^N$) becomes dense in the space
of all \K metrics as $r,N\to\infty$ \cite{Ti3}. Thus the quantum picture
tends to the classical one as $r\to\infty$.

\sec{The Yau-Tian-Donaldson conjecture.}
By analogy with the Kempf-Ness theorem in finite dimensions
(and by taking the infinite limit of Theorem \ref{KN}) it is natural to conjecture a \emph{Hitchin-Kobayashi correspondence} (the name coming from the analogy with the bundle case in the next section). That is a variety should admit a cscK metric if and only if it is polystable in a certain sense.

In fact Yau \cite{Y2} first suggested that there should be a relationship between stability and the existence of KE metrics.
Tian \cite{Ti1} proved this for surfaces, introduced his notion of K-stability, and, building on his work with Ding \cite{DT}, showed it was satisfied by K\"ahler-Einstein manifolds \cite{Ti2}. The
definition of K-stability was generalised to more singular test configurations by Donaldson \cite{Do4} who also showed that cscK implies K-semistability
\cite{Do3}. So it was thought that \textbf{K-polystability}, as defined above,
should be the right notion to be equivalent to cscK.

Recent explicit examples \cite{ACGT} in the extremal metrics case (where
there is a similar conjecture due to Sz\'ekelyhidi \cite{Sz}) suggest that this should be strengthened to \textbf{analytic K-polystability},
allowing more general analytic (instead of just algebraic) test configurations.
In particular one should allow the line bundle $L$ over the test configuration
to be an $\R$-line bundle: an $\R$-linear combination (by tensor product)
of $\C^*$-linearised line bundles. So the most likely
Yau-Tian-Donaldson conjecture as things stand at the end of 2005 is the following.

\begin{conj} \label{YTD}
$(X,L)$ is analytically K-polystable $\iff (X,L)$ admits a cscK metric. This is unique up to the holomorphic automorphisms of $(X,L)$.
\end{conj}

This would be the right higher dimensional generalisation of the
uniformisation theorem for Riemann surfaces.

There is very little progress on this conjecture in the $\Rightarrow$ direction
except for projective bundles \cite{BdB, Ho, RT2} and Donaldson's deep work
on toric surfaces \cite{Do4}. In the KE case there are sufficient conditions for existence given by Tian's
$\alpha$-invariant \cite{Ti4} and Nadel's multiplier ideal sheaf \cite{Na}, but no one has successfully related these to stability. Part of the problem, quite apart
from the analytical difficulties, is that we do not have a good intrinsic
understanding of stability for varieties -- i.e. no one has successfully analysed
the Hilbert-Mumford criterion for varieties.

Summarising the status of the whole theory for varieties, we have the infinite dimensional analogue of the balanced condition
for points in $\PP^1$ (i.e. cscK metrics) and part of the relationship
to stability, but not the algebro-geometric description
of stability. That is, the Hilbert-Mumford criterion, giving the analogue of the multiplicity $<n/2$ condition for points in $\PP^1$, is missing.

\begin{center} \input{dg1.pstex_t} \end{center}

\section{Moduli of bundles over $(X,L)$}

For holomorphic bundles $E$ over a polarised algebraic variety $(X,L)$
there is a very similar story which is more-or-less completely worked out.
Again there are subtleties due to
different notions of stability, but for bundles for which Gieseker and slope
stability coincide, for simplicity (such as those with coprime rank and degree,
or bundles over curves), we have, for $r\gg0$,

\begin{center} \input{dg2.pstex_t} \end{center}
\smallskip We now briefly explain this theory.

\sec{The gauge theory picture.}
The formal infinite dimensional picture was described by Atiyah-Bott
\cite{AB}. Fix a compatible hermitian metric on a $C^\infty$-bundle $E$ and consider the gauge
group $U(E)=\{$unitary $C^\infty$-maps $E\to E\}$ and its (genuine) complexification
$GL(E)$ of all $C^\infty$ invertible bundle maps $E\to E$. These act on
$$
\mathcal A=\big\{\text{unitary connections }A\text{ with }F_A^{0,2}=0\big\}. $$
The $U(E)$-action is obvious; $GL(E)$ acts by pulling back the $(0,1)$-part
$\dbar_A$ of the connection and then taking the unique Chern connection compatible
with both this and the metric. The integrability condition $F_A^{0,2}=\dbar_A^2=0$ ensures that $\dbar_A$ defines
a holomorphic structure on $E$. Thus any two $\dbar$-operators define isomorphic holomorphic structures on $E$ if and only if they lie in the same $GL(E)$-orbit.
So the formal
complex quotient of $\mathcal A$ by $GL(E)$ is the moduli space of holomorphic
vector bundles on $X$ with topological type $E$. (Of course we expect to need a stability condition to form this quotient.)

Alternatively, fixing the $\dbar$-operator and pulling back the metric by
$GL(E)$ gives the direct analogue of (\ref{corbit}) for the complexified
orbit of $\dbar_A$:
\beq{borbit}
U(E)\to GL(E).\dbar_A\twoheadrightarrow\big\{
\text{compatible metrics on }(E,\dbar_A)\big\}.
\eeq
The last map is onto since
$GL(E)$ acts transitively on the space of compatible hermitian metrics on
$E$ (the space of metrics being $GL(E)/U(E)$), so a complexified orbit can
be thought of as giving all compatible metrics on a fixed holomorphic bundle $(E,\dbar_A)$, up to the action of $U(E)$.

Fix a compatible hermitian metric
on $L$, inducing a \K form $\omega$ on $X$. Then $\mathcal A$ inherits a natural K\"ahler structure, with symplectic form given by $\Omega(a,b)=\int_X\tr(a\wedge
b)\wedge\omega^{n-1}$ for $a,b\in\Omega^1(\mathrm{End}\,E)$ tangent vectors
to $\mathcal A$. Atiyah-Bott show that $U(E)\acts\mathcal A$ has a moment
map
$$
A\mapsto F_A^{1,1}\wedge\omega^{n-1}-\lambda\id\omega^n\in\Omega^{2n}(\su(E)),
$$
thinking of the latter space as dual to $\Omega^0(\su(E))$ by the trace
pairing and integration. Here $\lambda=2\pi i\mu(E)\big/\int_X\omega^n$ is a topological constant, where
\beq{mu}
\mu(E)=\frac{\int_Xc_1(E).\omega^{n-1}}{\rk E}
\eeq
is the \emph{slope} of $E$.

Thus zeros of the moment map are \textbf{Hermitian-Yang-Mills} connections;
solutions of $\Lambda F_A^{1,1}=$\,const.id. An infinite dimensional version
of the Kempf-Ness theorem would be that in a polystable orbit of $GL(E)$
there should be a HYM connection (i.e. a metric whose associated Chern connection
is HYM; we call this a HYM metric), unique up to the action of $U(E)$, as conjectured by Hitchin and Kobayashi.

\begin{thm} \emph{[Donaldson-Uhlenbeck-Yau]} \label{DUY}
$E$ slope polystable $\Longleftrightarrow E$ admits a HYM metric. It is unique up to the automorphisms of $E$.
\end{thm}

The notion of stability that arises here (also called \emph{Mumford stability}) comes from GIT.

\sec{The GIT picture.}
Suppose we wanted form an algebraic moduli space of bundles $E$ over $(X,L)$
of fixed topological type. (More generally, to get a compact moduli space,
we have to consider coherent sheaves $E$ of the same Hilbert polynomial $\chi(E(r))$.)
We can twist $E(r):=E\otimes L^r,\ r\gg0$ until (a bounded subset of) the $E$s have no higher cohomology and are generated by their holomorphic sections:
\beq{quot1}
0\to\ker\to\underline{H^0(E(r))}\to E(r)\to0 \qquad\mathrm{on\ }X.
\eeq
Fixing an isomorphism $H^0(E(r))\cong\C^N,\ N=\chi(E(r))$, we have expressed
all such $E$s as quotients of $\O(-r)^{\oplus N}$.
Such quotients are easily parameterised algebraically by a Quot scheme (for
instance as a subset of the Grassmannian of subspaces $H^0(\ker(s))\le
H^0(\O(s))^{\oplus N}$ given by taking sections of (\ref{quot1}) tensored by $L^s,\ s\gg0$.)

So we divide by choices of the isomorphism $H^0(E(r))\cong\C^N$, i.e. by $SL(N,\C)$, to get the moduli space of (semistable) sheaves.

In this case the Hilbert-Mumford criterion can be manipulated (Mumford, Takemoto,
Gieseker, Maruyama, Simpson \cite{HL}) to give an algebro-geometric understanding of stability. We describe this in a later section; the upshot is the following.

We write the Hilbert polynomial $\chi(E(r))=a_0r^n+a_1r^{n-1}+\ldots$, where
$a_0=\rk E\int_X\omega^n/n!,\ \ a_1=\int_Xc_1(E).\omega^{n-1}/(n-1)!
+\varepsilon(X)$, etc. are topological.
We use its monic version, the reduced Hilbert polynomial
\beq{rH}
p\_E(r):=\frac{\chi(E(r))}{a_0}=r^n+\frac{a_1}{a_0}r^{n-1}+\ldots\,.
\eeq

Then $E$ is stable if and only if \emph{for all} coherent subsheaves
$F\into E,\ p\_F(r)\prec p\_E(r)$ in the following sense (depending on the line bundle chosen on the Quot scheme):
\begin{itemize}
\item \emph{Gieseker stable} if and only if $p\_F(r)<p\_E(r)\quad\forall
r\gg0$.
\item \emph{Slope stable} if and only if $\frac{a_1(F)}{a_0(F)}<
\frac{a_1(E)}{a_0(E)}\ \big(\!\iff\mu(F)<\mu(E)\big)$.
\end{itemize}
Here, as before, $\mu(E)=\int_Xc_1(E).\omega^{n-1}/\rk(E)$ is the slope of
$E$ (\ref{mu}). Gieseker and slope stability coincide on curves $X$. Slope stability corresponds to taking a certain semi-ample line bundle on the Quot scheme (roughly speaking given by restricting sheaves to high degree complete intersection curves in $X$ and using the usual line bundle for moduli of bundles on the curve). GIT needs
amending for this situation; so far this has been carried out only for $X$ a surface by Jun Li \cite{HL}.

Semistability is similar (replacing $<$ by $\le$), while polystability is
equivalent to semistability where the only destabilising subsheaves $F$ are
direct summands of $E$. Slope polystability then turns out to be the right
stability notion for the infinite dimensional quotient and HYM of Theorem \ref{DUY}.

\sec{The symplectic reduction picture.}
For $E$ a bundle, we can interpret the quotient $\underline{\C}^N\to E(r)\to0$ (\ref{quot1}) differently, via its classifying map $X\to Gr$ to the Grassmannian of quotients of $\C^N$.

Then, much as in the varieties case, we fix compatible hermitian metrics
on $L$ and $E$ inducing a \K form on $X$ and an $L^2$-metric on $\C^N\cong
H^0(E(r))$.
Then there are actions of $SU(N)< SL(N,\C)$ on $Gr$, inducing
a moment map $m\colon Gr\into\su(N)^*$ and an action $SU(N)\acts$ Maps($X,Gr)$.
Its moment map is the integral of (the pullback of) $m$ over $X$, so we can again talk about \emph{balanced} $X\to Gr$ (those with centre of mass zero
in $\su(N)^*$) and asymptotics as $r,N\to\infty$.

Proving conjectures of Donaldson \cite{Do2}, Wang shows that the existence of a balanced map is equivalent to the Gieseker
polystability of $E$ \cite{Wa1}. Slope stable bundles (which are therefore
Gieseker stable) admit balanced maps $X\to Gr$ for $r\gg0$ \cite{Wa2}, and pulling back the canonical quotient connection on $Gr$ and taking $\lim_{r\to\infty}$ gives a conformally Hermitian-Yang-Mills connection on $E$ (which is HYM
after rescaling). (Unfortunately, this is not how the results are proved;
Wang uses the Donaldson-Uhlenbeck-Yau theorem to give an a priori HYM connection
which can be compared to the sequence of balanced metrics.)

\section{Slope criteria for algebraic varieties}

\sec{Slope for K-stability.}
So we have seen that the finite and infinite dimensional GIT and symplectic reduction pictures work for bundles, and tend to one another as $r\to\infty$. In particular the stability notion for bundles and sheaves involves only subsheaves $F< E$. So one might ask if the subvariety $\PP(F)\subset\PP(E)$ can destabilise
the variety $\PP(E)$. Or if, more generally, subschemes $Z\subset(X,L)$ can
destabilise $(X,L)$. (Cf. length $\ge n/2$ subschemes destabilising length-$n$ schemes in $\PP^1$ (\ref{P1}).)

We need some topological data for $(X,L)$ analogous to that for bundles (\ref{rH}).
Fixing $Z\subset(X,L)$, we have the Hilbert polynomial of $L=\O_X(1)$
$$
h^0(\O_X(r))=a_0r^n+a_1r^{n-1}+\ldots\,,
$$
and the \emph{Hilbert-Samuel} polynomial of $\I_Z$ (for $x\in\Q$ and
$rx\in\mathbb N$):
$$
h^0(\I_Z^{xr}(r))=a_0(x)r^n+a_1(x)r^{n-1}+\ldots\,.
$$
Then by working on the blow up of $X$ in $Z$ one can see that the $a_i(x)$
are polynomials in $x\in\Q\cap[0,\epsilon(Z))$ for $r\gg0$. (More precisely,
there is a constant $p$ and a polynomial which is equal to $a_i(x)$ for $xr>p$
or $x=0$.) Here $\epsilon(Z)$ is the \emph{Seshadri constant} of $Z$, the
supremum of $x$ such that $\I_Z^{xr}(r)$ is generated by global sections for $r\gg0$.

For $X$ normal, $a_0(0)=a_0$, and $a_1(0)=a_1$. All the $a_i(x)$ are given
by topological formulae by Riemann-Roch, for instance
$$
a_0=\frac{\int_X\omega^n}{n!}\,, \qquad
a_1=\frac{\int_Xc_1(X)\omega^{n-1}}{2(n-1)!}\,.
$$

For any $c\le\epsilon(Z)$, analogously to the definition of slope for bundles
(\ref{mu}), we define the slope of $Z$ to be
\beq{kslope}
\mu_c(\I_Z)=\frac{\int_0^ca_1(x)+\frac{a_0'(x)}2dx}{\int_0^ca_0(x)dx}\,.
\eeq
$Z=\emptyset$ gives
$$
\mu(X)=\frac{a_1}{a_0}\,.
$$
We then have the following \cite{RT2}.

\begin{thm} \label{kkss}
$(X,L)$ K-semistable $\Rightarrow$ slope semistable: $\mu_c(\I_Z)\le\mu(X)$ for all closed subschemes $Z\subset X$ and $c\le\epsilon(Z)$.
\end{thm}

Removing the ``semi" is a little more involved. If we use the algebraic
K-stability of \cite{Do4, RT1} then we define slope stability to mean
\begin{itemize}
\item $\mu_c(\I_Z)<\mu(X)\ \ \forall c\in(0,\epsilon(Z))\cap\Q$, \emph{and}
\item $\mu_{\epsilon(Z)}(\I_Z)<\mu(X)$ if $\epsilon(Z)\in\Q$ and
$\I_Z^{\epsilon(Z)r}(r)$ is saturated by global sections for $r\gg0$.
\end{itemize}
(The quickest definition of saturated \cite{RT1} is that on the blow up $\pi\colon\Bl_ZX\to X$ of $X$ in $Z$ with exceptional divisor $E$, $\pi^*\I_Z^{\epsilon(Z)r}(r)=\pi^*L^r(-\epsilon(Z)r)$ should be generated by global sections.) Similarly slope polystability is defined as slope stability
except in the second part of the definition we allow $\mu_{\epsilon(Z)}(\I_Z)$ to equal $\mu(X)$ if on $\Bl_{Z\times\{0\}}(X\times\C),\ 
L(-\epsilon(Z)P)$ (where $P$ is the exceptional divisor) is pulled back from a contraction $\Bl_{Z\times\{0\}}(X\times\C)\to X\times\C$ (which of course won't be the original blowup map).

If we use analytic K-stability \cite{RT2}, which is what should be relevant
to the cscK problem (Conjecture \ref{YTD}), then the relevant definition of slope stability allows irrational $c$\,:
\begin{itemize}
\item $\mu_c(\I_Z)<\mu(X)\ \ \forall c\in(0,\epsilon(Z))$, \emph{and}
\item $\mu_{\epsilon(Z)}(\I_Z)<\mu(X)$ if $\epsilon(Z)\in\Q$ and
$\I_Z^{\epsilon(Z)r}(r)$ is saturated by global sections for $r\gg0$.
\end{itemize}
Again slope polystability is defined in the same way except for the second
condition in which we allow $\mu_{\epsilon(Z)}(\I_Z)=\mu(X)$ if on
$\Bl_{Z\times\{0\}}(X\times\C),\ 
L^r(-\epsilon(Z)rP)$ is pulled back from a contraction to $X\times\C$.

Then the analogue of Theorem \ref{kkss}, for the either notion of K-stability,
is the following \cite{RT1}.

\begin{thm} \label{thm2}
$(X,L)$ K-(poly)stable $\Rightarrow$ slope (poly)stable.
\end{thm}

As a corollary of Theorem \ref{kkss} and the results of Donaldson and Chen-Tian
mentioned in Theorem \ref{thm} we find the following.

\begin{cor}
If $\mu_c(\I_Z)>\mu(X)$ then $X$ does not admit a cscK metric in the class
of $c_1(L)$.
\end{cor}

We give some examples.
\begin{itemize}
\item If $F<E$ is a slope-destabilising subbundle of a vector bundle $E\to X$ the $\PP(F)\subset\PP(E)$ destabilises for suitable polarisations $\pi^*L^m\otimes\O_{\PP(E)}(1),\
m\gg0$, on $\PP(E)$. Conversely \cite{Ho}, if $E$ is slope-stable and the
base $X$ is cscK with discrete automorphism group then $\PP(E)$ is cscK for $m\gg0$.
\item When the base is a curve, we can do better \cite{RT2}. In fact, for $\PP(E)$ with discrete automorphism group and \emph{any} polarisation,
$$
\PP(E) \ \mathrm{cscK} \iff E\ \mathrm{HYM} \iff E\ \mathrm{stable}.
$$
The converses follow from the Narasimhan-Seshadri theorem \cite{BdB} that stable bundles admit projectively flat connections.
\item $-1$-curves destabilise del Pezzo surfaces $(X,L)$ for appropriate $L$. In particular one can find examples with reductive (even trivial) automorphism
group, showing that the folklore conjecture Aut\,$(X)$
reductive $\Rightarrow$ cscK does not hold for surfaces. (Tian showed that
it does not hold for threefolds \cite{Ti2}, but that it does hold for
anticanonically polarised surfaces \cite{Ti1}.)
\item For instance $\PP^2$ blown up in one point cannot admit a cscK metric
for any polarisation because Aut$\,(X)$ is not reductive. From the above
point of view this is because it is destabilised by the exceptional $-1$-curve for all polarisations.
\item Generically stable varieties can specialise to unstable ones. For instance
moving two $-1$-curves together on a stable del Pezzo gives a limiting $-2$-curve
(if the $-1$-curves arise from blowing up distinct points, then blow up two
``infinitely near" points -- one point and then another on the exceptional
divisor) which can destabilise for suitable $L$.
\item Calabi-Yau manifolds, and varieties with canonical singularities
and numerically trivial canonical bundle ($mK_X\sim0$) are slope stable.
\item Canonically polarised varieties with canonical singularities
(i.e. the canonical models of Mori theory) are slope stable.
\item Remarkably, the product of a (nongeneric) smooth curve with itself can admit polarisations which are slope unstable, giving surfaces of general
type which are neither K- nor Chow stable, and which do not admit cscK metrics in that class \cite{Ro}.
\end{itemize}

We will sketch the algebro-geometric proof of the above results later. But
differential-geometrically, what the proofs amount to is the following.

We know that the Mabuchi energy (or log-norm functional, in GIT language) is convex over the space of all K\"ahler metrics on our fixed complex manifold
$X$. Intuitively, if it is proper in some sense (roughly, tends to $+\infty$
at infinity) then it should have a unique absolute minimum, which, modulo
regularity issues, would be cscK. Conversely the manifold is strictly K-unstable,
with no cscK metric, if the Mabuchi energy is unbounded below. If $X$ is
slope destabilised by some subscheme $Z$ then there is a family of K\"ahler
metrics on $X$, given by ``stretching the neck" around $Z$, along which the
Mabuchi energy tends to $-\infty$, so $X$ cannot be cscK.

\sec{Slope for Chow stability}. Fix $Z\subset(X,\O_X(1))\subset(\PP^N,\O(1))$ embedded by sections of $\O_X(1)$, and, as before,
$$
h^0(\O_X(r))=a_0r^n+a_1r^{n-1}+\ldots\,,
$$$$
h^0(\I_Z^{xr}(r))=a_0(x)r^n+a_1(x)r^{n-1}+\ldots\,.
$$
Then, for all $c\le\epsilon(Z)\in\mathbb N$, define the \emph{Chow slope} of $Z$ to be
$$
Ch_c(\I_Z)=\frac{\sum_{i=1}^ch^0(\I_Z^i(1))}{\int_0^ca_0(x)dx}\,,
$$
a discrete version of (\ref{kslope}). $Z=\emptyset$ gives
$$
Ch(X)=\frac{h^0(\O_X(1))}{a_0}=\frac{N+1}{a_0}\,,
$$
and we have the following \cite{RT1}.

\begin{thm} \label{css}
Chow (semi)stable $\Rightarrow$ slope (semi)stable: $$Ch_c(\I_Z)\mathop{<}_{(\le)}
Ch(X)\ \ \forall Z\subset X.$$
\end{thm}

To see where these results come from, and to explain how far one can get
towards a converse, we need to analyse the Hilbert-Mumford criterion for
varieties $(X,L)$. We first warm up with a brief overview of the bundle case.

\sec{The Hilbert-Mumford criterion for vector bundles.}
Given a coherent sheaf $E$ on $X$, recall (\ref{quot1}) how an picking an identification $H^0(E(r))\cong\C^N$ for $r\gg0,\,N=\chi(E(r))$, realises $E(r)$ as a point of a Quot scheme of quotients
\beq{quot}
\O(-r)^{\oplus N}\to E\to0.
\eeq
Dividing out by the identification, i.e. quotienting the relevant subset of Quot by $SL(N,\C)$, gives a moduli space of sheaves.

To apply the Hilbert-Mumford criterion, we consider the $E$-orbit of a 1-PS $\C^*< GL(N,\C)$ on Quot \cite{HL} (we shall normalise to $SL(N,\C)$ later), whose eigenvalues we can assume are all positive, without loss of
generality. The eigenspaces $V_\lambda<\C^N\cong H^0(E(r))$
give a weight filtration $V_{\le i}=\oplus_{\lambda\le i}V_\lambda$ of $\C^N$.
Their images $F_i\le E$ under the map (\ref{quot}) give a filtration of $E$,
$$
F_0\le F_1\le \ldots\le F_p\le E,
$$
and the orbit gives a sheaf over $X\times\C$ described in terms of the $F_i$
as follows. Let $\mathbb F_i,\,\mathbb E$ denote the pullbacks of the sheaves
$F_i,\,E$ to $X\times\C$, so $\mathbb F_i=F_i\otimes\C[t]$ where $t$ is the variable pulled back from $\C$. Then the orbit gives the following subsheaf
of $\mathbb E$,
\beq{bideal}
\mathbb F_0+t.\mathbb F_1+t^2.\mathbb F_2+\ldots+t^p.\mathbb F_p+t^{p+1}.\mathbb E< E\otimes\C[t].
\eeq
This is a degeneration of the general fibre $E$ to
\beq{decom}
F_0\ \oplus\ F_1/F_0\ \oplus\ \ldots\ \oplus\ F_p/F_{p-1}\ \oplus\ E/F_p
\eeq
over the central fibre $X\times\{0\}$. One can prove this inductively as
follows.

Set $Q_i=E/F_i$, giving exact sequences $0\to F_i/F_{i-1}\to Q_{i-1}\to
Q_i\to0$. In
the $p=0$ case, $\mathbb E_0:=\mathbb F_0+t.\mathbb E$ is the kernel of the composition $\mathbb E\to E\to Q_0$, where the latter two sheaves are considered to be supported on the central fibre $X\times\{0\}$. (So $\mathbb E_0$
is the \emph{elementary transform} of $\mathbb E$ in $Q_0$ on $X\times\{0\}$.) $\mathbb E_0$ has a map to $F_0$ induced by the diagram
\beq{elem}
\xymatrix{
0\rto & \mathbb E_0\rto\dto & \mathbb E\rto\dto & Q_0\rto\ar@{=}[d] & 0 \\
0\rto & F_0 \rto & E\rto & Q_0\rto
& 0\,,\!\!}
\eeq
and one to $Q_0$:
$$
\xymatrix{
0\rto & \mathbb E_0\rto\dto & \mathbb E\rto\dto & Q_0\rto\ar@{=}[d] & 0 \\
0\rto & Q_0 \rto^(0.35){.t} & Q_0\otimes\frac{\C[t]}{(t^2)}\rto & Q_0\rto
& 0\,.\!\!}
$$
The pair make the central fibre of $\mathbb E_0$ isomorphic to $F_0\oplus
Q_0$, while the map $\mathbb E_0\to Q_0$ can be composed with the surjection
$Q_0\to Q_1$ to continue the induction by defining $\mathbb E_1$ as its kernel.
This is clearly just $\mathbb F_0+t.\mathbb F_1+t^2.\mathbb E$, and similar
working shows it has central fibre $F_0\oplus F_1/F_0\oplus Q_1$, and so on.

Different 1-PSs can give the same filtration (if the $i$th piece of the weight
filtration of $\C^N$ generates the same subsheaf $F_i\le E$). But for every filtration
$F_i$ and sequence of weights there is a canonical \emph{least stable} 1-PS, given by choosing the weight filtration to be $V_{\le i}=H^0(F_i(r))< H^0(E(r))$.

So we need only consider these canonical 1-PSs: they have the largest GIT
weight in their class. Since (\ref{decom}) gives
the weight space decomposition of the limiting sheaf over the central fibre, the weights of these 1-PSs $\sum_nnh^0\big((F_n/F_{n-1})(k)\big),\ k\gg0,$ are positive linear combinations of weights of the canonical 1-PSs associated to the splittings
\beq{EF}
F_i\ \oplus\ E/F_i.
\eeq
So in fact we need only control the weights of these simpler splittings.
Calculating their weights (and then normalising into $SL(N,\C)$ \cite{HL}) gives the Gieseker stability condition for bundles described earlier, controlled by single subsheaves $F<E$ and their reduced Hilbert polynomials.

\sec{The Hilbert-Mumford criterion for varieties.} We can now try to do the analogous thing for varieties, following \cite{Mu, RT1}.

Any test configuration $(\X,\L)$ is $\C^*$-birational to $(X\times\C,L)$,
so is (a contraction $p$ of) the blow up of $X\times\C$ in a $\C^*$-invariant
ideal $I$ supported on the central fibre, with polarisation $p^*\L=
(\pi^*L)(-cP)=\pi^*(L\otimes I^c)$, where $P$ is the exceptional divisor.
\beq{testc}
\begin{array}{ccccc}
&& \!\!\Bl_I(X\times\C)\!\! \\
& {}^\pi\!\!\swarrow && \searrow^p \\
X\times\C \!\!\!\! &&&& \!\!\!\!\! \X
\end{array}
\eeq
Classifying such $\C^*$-invariant ideals, there exist subschemes
\beq{P}
Z_{p-1}\subseteq\ldots\subseteq Z_1\subseteq
Z_0\subseteq X
\eeq
with ideal sheaves $\I_{p-1}\supseteq\ldots\supseteq\I_1
\supseteq\I_0$ such that \cite{Mu}
\beq{ideal}
I=\I_0+t\I_1+t^2\I_2+\ldots+t^{p-1}\I_{p-1}+t^p.
\eeq

\bigskip \begin{flushleft}
\input{smallsubscheme.pstex_t}
\end{flushleft} \bigskip

Firstly we have an analogue of the fact that, in the bundle case, one need
only consider canonical 1-PSs. Namely, under any map of test configurations
like $p$ above, the weights are less stable (more positive) on the dominating
test configuration. (Notice that the blow up map $\Bl_I(X\times\C)\to X\times\C$
above is not \emph{such} a map of test configurations, since it does not
preserve polarisations: the line bundle is not a pullback from downstairs.)

\begin{prop} \label{stein}
Suppose $X$ is normal. Given a test configuration $(\X,\L)$ for $(X,L)$ and another flat $\C^\times$-family $\Y\to\C$ with a birational $\C^*$-equivariant map $p\colon\Y\to\X$, there exists an $a\ge0$ such that
$$
w(H^0_{\X_0}\!(\L^k))
=w\left(H^0_\Y(p^*\L^k)\big/tH^0_\Y(p^*\L^k)\right)-ak^n+O(k^{n-1}).
$$
\end{prop}

(Here $w$ denotes the total weight of a $\C^*$-action -- i.e. its weight on the determinant of the $\C^*$-module. The normality of $X$ is required to equate $H^0(X,L^k)$ with sections of $p^*\L^k$ on a general fibre of $\Y$.
The result is stated in rather more generality than we require here (in particular
allowing $\Y$ to have general fibre some blow up of $X$, rather than just
$X$) for future use. We are forced to use $H^0_\Y(p^*\L^k)\big/tH^0_\Y(p^*\L^k)$, rather than $H^0_{\Y\_0}(p^*\L^k)$, because $p^*\L$ need not be ample on $\Y$.)

So we need only consider weights on the normalisation of the blow up of $X\times\C$
in $I$, as this is is itself a perfectly good test configuration. (In the
Chow stability case this test configuration may not arise from a linear transformation
of the given projective space -- only from an embedding by a higher twist
of $L$ -- but this does not concern us as Proposition \ref{stein} gives an equality in weights to $O(k^n)$, which is all that is required for Chow stability.)

Next we consider the easiest case of $p=1$, i.e. $I=\I_0+(t)$. So we blow up in $Z_0\times\{0\}\subset X\times\C$, giving the
\emph{deformation to normal cone of} $Z_0$.
This is the analogue of the simplest degenerations of bundles and sheaves earlier; the canonical degenerations of a one-step filtration to a direct
sum (\ref{EF}).

\begin{center}
\input{smallnormalcone0.pstex_t}
\end{center}

The exceptional divisor $P$ is the \emph{normal cone} of $Z_0$: if $Z_0\subset X$ is smooth then this is the projective completion
$\PP(\nu\_{Z_0}\oplus\underline\C)\to Z_0$ of the normal bundle $\nu\_{Z_0}\to
Z_0$.

$\C^*\ni\lambda$ acts on the blow up (as $[1:\lambda]=[\lambda^{-1}:1]$ on
$\PP(\nu\_{Z_0}\oplus\underline\C)$ in the
smooth case) and on the line bundle $\pi^*L(-cP)$ over $\Bl_{Z_0\times\{0\}}
(X\times\C)$.

This deformation to the normal cone of $Z$ replaces $H^0_X(L^r)$ (filtered by $H^0(L^r\otimes\I^j_Z)$) by the associated graded
of the filtration on the central fibre:
\begin{multline*}
H^0_X(\I_Z^{cr}(r))\oplus tH^0_X\left(\I_Z^{cr-1}(r)\big/
\I_Z^{cr}(r)\right)\oplus\ldots \\
\oplus\ t^{cr-1}H^0_X\left(\I_Z(r)\big/
\I_Z^2(r)\right)\oplus t^{cr}H^0_X(\O_Z(r)).
\end{multline*}
Here $t$ is the coordinate on the $\C$-base. This is the splitting of sections of $L^r(-crP)$ on the central fibre
into those on the proper transform of the original central fibre $X$ (the
first term) plus the polynomials on $P$. In turn the latter can be split
into their $\C^*$-weight spaces as $t^j$ times by the homogeneous polynomials on the normal bundle of $Z\subset X$ of degree $cr-j$. 

So this is the weight space decomposition, with $\C^*$ acting on $t^j$ with weight $-j$, and the weight on top exterior power is
\begin{eqnarray}
w_r&=&-\sum_{j=1}^{cr}jh^0_X\left(\I_Z^{cr-j}(r)\big/
\I_Z^{cr-j+1}(r)\right) \nonumber\\
&=&-\sum_{j=1}^{cr}h^0(\I_Z^j(r))-crh^0(\O_X(r)). \label{cw}
\end{eqnarray}
This looks like a discrete approximation (Riemann sum) for $\int_0^ch^0(\I_Z^{xr}(r))dx$, so we estimate it by the trapezium rule, giving
$$
-\!\left(\int_0^c\!\!a_0(x)dx\!\right)\!r^{n+1} + \int_0^c\!
\left(\!a_1(x) + \frac{a_0'(x)}{2}\right)\!dx\ r^n+O(r^{n-1}).
$$
Normalising (to make the 1-PS lie in $SL(N,\C)$ instead of $GL(N,\C)$) we
find the K-stability slope criterion of (\ref{kslope}) and Theorems \ref{kkss} and \ref{thm2}.

Alternatively, taking the leading order term of (\ref{cw}) and normalising
gives the Chow slope criterion of Theorem \ref{css}.

\sec{General test configurations.} So we have the analogue of the result
in the bundle case that one need only consider canonical 1-PS orbits -- i.e.
we need only consider test configurations that are normalisations of the form (\ref{testc}). And we have the analogue of the simple degenerations
(\ref{EF}), given by the deformation to the normal cone of a subscheme $Z\subset
X$ yielding the right analogue of slope stability. For bundles any canonical 1-PS weight turned out to be a positive linear combination of these simple
weights, so ideally one would like to write the weight of a test configuration (\ref{testc}) as a positive linear combination of weights (\ref{cw}) -- if so one could conclude that stability of varieties was equivalent to slope stability.

So we need the analogue of the induction (\ref{elem}) that we did in the bundle case to handle
(\ref{bideal}). In fact it turns out we can mirror it almost completely;
moreover the bundle induction is a special case of what follows below when
we consider the sheaf $\mathbb E$ of (\ref{bideal}) to be the sections
of the polarisation $\O_{\PP(E^*)}(1)$ on the variety $\PP(E^*)\times\C$. The correspondence between (\ref{bideal}) and (\ref{ideal}) is clear, and the elementary transformations we did in the bundle case become, on projectivisation, the blowups below.

The proper transform $\overline{Z_0\times\C}$ of $Z_0\times\C$ is \emph{flat}
over the base $\C$. It defines a copy $Z_0'$ of $Z_0$ in the central fibre
of the deformation to the normal cone of $Z_0$:

\begin{center}
\input{normalcone.pstex_t}
\end{center}
 
By (\ref{P}) this defines subschemes $Z_{p-1}'\subseteq\ldots\subseteq Z_1'\subset Z_0'$. So next we blow up in $Z_1'$, giving $Z_{p-1}''\subseteq\ldots\subseteq Z_1''$; next blow up $Z_2''$, and so on inductively up to $Z_{p-1}^{(p-1)}$.

(For $X=\PP(E^*)$ and the degeneration (\ref{bideal}), we are blowing up the subschemes $\PP(Q_0^*)$ (in the central fibre), then the central fibre of the proper transform of $\PP(Q_1^*)$, and so on; but this is just the projectivisation of the elementary transformations (\ref{elem}).)

\begin{thm} The blow up of $X\times\C$ in $I=\I_0+t\I_1+\ldots+t^{p-1}\I_{p-1}+t^p$
is a contraction of this iterated blow up.
\end{thm}

This is meant in the $\C^*$-equivariant polarised sense -- the polarisation $L(-cP)$ on $\Bl_I(X\times\C)$
is the pullback of the natural polarisation on the iterated blow up given
by starting with $L$ and, at each stage, pulling back and subtracting $c$ times the exceptional divisor.

Thus, by our result (Proposition \ref{stein}) that one need only calculate weights on dominating test
configurations, we are left with calculating the amount that each blow up
in $Z_i^{(i)}$ adds to the weight of the resulting $\C^*$-action on the determinant
of the space of sections of the $r$th power of the polarisation on the central
fibre.

\begin{thm} \label{train}
Consider the $i$th step, when we blow up $Z_i^{(i)}$. If all thickenings of $\overline{(Z_i\times\C)}$ are flat over $\C$ then this adds $w(Z_i)$ to the weight, to $O(r^n)$. \\
(Here $w(Z_i)$ is weight on deformation to normal cone of $Z_i$.)
\end{thm}

In fact under certain conditions one can get the result to $O(r^{n-1})$ \cite{RT1}.

So if this flatness condition holds, the total (normalised) weight is $w(Z_0)+\ldots+w(Z_{p-1})$.
So $X$ is stable if and only if
\begin{equation} \label{result}
w(Z_0)+\ldots+w(Z_{p-1})\prec0\iff w(Z)\prec0\ \forall Z.
\end{equation}
So if this flatness condition held in general, then stability and slope stability would be equivalent.

\sec{The flatness problem.}
In fact $\overline{(Z_i\times\C)}$ \emph{is} flat over $\C$, but since blow ups use all powers of an ideal, we require all of its scheme theoretic thickenings $k\overline{(Z_i\times\C)}$
(defined by the ideal $\I_{\overline{Z_i\times\C}}^k$) to be flat too.
The idea of the proof of Theorem \ref{train} is then that a formal neighbourhood of $Z_i^{(i)}$ in the test configuration looks sufficiently like a formal neighbourhood of $Z_i\times\{0\}\subset
X\times\C$ for the weights added by the two blow ups to be comparable to
$O(r^n)$ -- the corrections being due to estimates on the sizes of nonvanishing
cohomology groups.

This flatness condition holds for $Z_i\subset X$ smooth, or reduced simple normal crossing (snc) divisors.

In general, one can use resolution of singularities
$$
(X\supset Z_i)\stackrel{\,\pi}\longleftarrow(\widehat X\supset m_iD_i),\ \ D_i\ \mathrm{snc\ divisors},
$$
to replace $(X,L)$ by $(\widehat X,\pi^*L)$. Test configurations for the
latter dominate those of the former, so Proposition \ref{stein} allows us
to obtain (\ref{result}) for $X$ \emph{normal}, so long as $m_i=1$ for all $i$.

So finally we need to be able to deal with the snc divisors $D_i$ being possibly
nonreduced. We can attempt to do this by basechange \cite{RT1}, which we illustrate with an example.

Consider the case where $Z_0\subset X$ is a double point in a smooth curve. Locally then $\I_0=(x^2)$, and the deformation to normal cone of $Z_0$ is
the blow up of $X\times\C$ in $(x^2,t)$.

Now consider squaring the $\C^*$-action. This is trivial from
a GIT point of view (it just doubles the weight, which we can
later halve). But it fundamentally alters our geometric description of the
test configuration, blowing up $X\times\C$ in the ideal $(x^2,t^2)$.

Taking the integral closure of this ideal corresponds to normalising the blow up, which we can deal with by Proposition \ref{stein}. That is, we get a more unstable
test configuration by blowing up in $(x^2,xt,t^2)=(x,t)^2$, and it suffices
to control the weights of this test configuration.

But this is now a much nicer ideal, and corresponds to blowing up in $(x,t)$ and using a different line bundle (squaring the exceptional divisor $P\mapsto2P$ or $c\mapsto2c$). So modulo doubling $c$ and the weight, we have removed
the multiplicity 2 of the double point $Z_0$.

In this way we can deal with $D_i$ with multiplicities $m_i$ when they
\emph{all have the same support}. This is enough to prove that (K- and Chow) stability coincides with (K- and Chow) slope stability for smooth curves, and indeed gives probably the ``right" geometric proof, rather than the old combinatorial (for Chow stability) and analytical (for K-stability) proofs.

Of course, for higher dimensions, one would like to combine the two approaches to deal with snc divisors with intersecting components of different multiplicities, for instance, $D_0=(x^2y=0),\ D_1=(x=0)$. This is still work in progress,
but its difficulty suggests that  slope stability is not enough to describe
stability (unlike for
the more linear bundle case) except in one dimension or on projective bundles
over stable bases.

\end{document}

%% file: blowup.pstex_t
\begin{picture}(0,0)%
\includegraphics{blowup.pstex}%
\end{picture}%
\setlength{\unitlength}{2368sp}%
\begingroup\makeatletter\ifx\SetFigFont\undefined%
\gdef\SetFigFont#1#2#3#4#5{%
  \reset@font\fontsize{#1}{#2pt}%
  \fontfamily{#3}\fontseries{#4}\fontshape{#5}%
  \selectfont}%
\fi\endgroup%
\begin{picture}(4662,3729)(2389,-6448)
\put(4576,-2911){\makebox(0,0)[lb]{\smash{{\SetFigFont{11}{13.2}{\rmdefault}{\mddefault}{\updefault}{\color[rgb]{0,0,0}$X$}%
}}}}
\put(6901,-3661){\makebox(0,0)[lb]{\smash{{\SetFigFont{11}{13.2}{\rmdefault}{\mddefault}{\updefault}{\color[rgb]{0,0,0}$\O_X(-1)$}%
}}}}
\put(7051,-5836){\makebox(0,0)[lb]{\smash{{\SetFigFont{11}{13.2}{\rmdefault}{\mddefault}{\updefault}{\color[rgb]{0,0,0}$\widetilde X$}%
}}}}
\end{picture}%

%% file: orbits.pstex_t
\begin{picture}(0,0)%
\includegraphics{orbits.pstex}%
\end{picture}%
\setlength{\unitlength}{4144sp}%
\begingroup\makeatletter\ifx\SetFigFont\undefined%
\gdef\SetFigFont#1#2#3#4#5{%
  \reset@font\fontsize{#1}{#2pt}%
  \fontfamily{#3}\fontseries{#4}\fontshape{#5}%
  \selectfont}%
\fi\endgroup%
\begin{picture}(886,1328)(1924,-1782)
\end{picture}

%% file: smallHM.pstex_t
\begin{picture}(0,0)%
\includegraphics{smallHM.pstex}%
\end{picture}%
\setlength{\unitlength}{3315sp}%
\begingroup\makeatletter\ifx\SetFigFont\undefined%
\gdef\SetFigFont#1#2#3#4#5{%
  \reset@font\fontsize{#1}{#2pt}%
  \fontfamily{#3}\fontseries{#4}\fontshape{#5}%
  \selectfont}%
\fi\endgroup%
\begin{picture}(6304,2260)(1711,-2198)
\put(4738,-1905){\makebox(0,0)[lb]{\smash{{\SetFigFont{11}{13.2}{\rmdefault}{\mddefault}{\updefault}{\color[rgb]{0,0,0}$X$}%
}}}}
\put(2189,-2144){\makebox(0,0)[lb]{\smash{{\SetFigFont{11}{13.2}{\rmdefault}{\mddefault}{\updefault}{\color[rgb]{0,0,0}$\C^*.x$}%
}}}}
\put(3511,-691){\makebox(0,0)[lb]{\smash{{\SetFigFont{11}{13.2}{\rmdefault}{\mddefault}{\updefault}{\color[rgb]{0,0,0}$\tilde x$}%
}}}}
\put(2926,-106){\makebox(0,0)[lb]{\smash{{\SetFigFont{11}{13.2}{\rmdefault}{\mddefault}{\updefault}{\color[rgb]{0,0,0}$\O_X(-1)$}%
}}}}
\put(1711,-691){\makebox(0,0)[lb]{\smash{{\SetFigFont{11}{13.2}{\rmdefault}{\mddefault}{\updefault}{\color[rgb]{0,0,0}$\C^*.\tilde x$}%
}}}}
\put(4231,-2086){\makebox(0,0)[lb]{\smash{{\SetFigFont{11}{13.2}{\rmdefault}{\mddefault}{\updefault}{\color[rgb]{0,0,0}$x_0$}%
}}}}
\put(3556,-1996){\makebox(0,0)[lb]{\smash{{\SetFigFont{11}{13.2}{\rmdefault}{\mddefault}{\updefault}{\color[rgb]{0,0,0}$x$}%
}}}}
\put(4726,-961){\makebox(0,0)[lb]{\smash{{\SetFigFont{10}{12.0}{\rmdefault}{\mddefault}{\updefault}{\color[rgb]{0,0,0}$\O_{x_0}(-1)$}%
}}}}
\put(7246,-511){\makebox(0,0)[lb]{\smash{{\SetFigFont{12}{14.4}{\rmdefault}{\mddefault}{\updefault}{\color[rgb]{0,0,0}stable}%
}}}}
\put(7246,-1006){\makebox(0,0)[lb]{\smash{{\SetFigFont{12}{14.4}{\rmdefault}{\mddefault}{\updefault}{\color[rgb]{0,0,0}semistable}%
}}}}
\put(7246,-1501){\makebox(0,0)[lb]{\smash{{\SetFigFont{12}{14.4}{\rmdefault}{\mddefault}{\updefault}{\color[rgb]{0,0,0}unstable}%
}}}}
\end{picture}%

%% file: smallorb.pstex_t
\begin{picture}(0,0)%
\includegraphics{smallorb.pstex}%
\end{picture}%
\setlength{\unitlength}{3315sp}%
\begingroup\makeatletter\ifx\SetFigFont\undefined%
\gdef\SetFigFont#1#2#3#4#5{%
  \reset@font\fontsize{#1}{#2pt}%
  \fontfamily{#3}\fontseries{#4}\fontshape{#5}%
  \selectfont}%
\fi\endgroup%
\begin{picture}(3477,2691)(1711,-2198)
\put(4738,-1905){\makebox(0,0)[lb]{\smash{{\SetFigFont{11}{13.2}{\rmdefault}{\mddefault}{\updefault}{\color[rgb]{0,0,0}$X$}%
}}}}
\put(2189,-2144){\makebox(0,0)[lb]{\smash{{\SetFigFont{11}{13.2}{\rmdefault}{\mddefault}{\updefault}{\color[rgb]{0,0,0}$\C^*.x$}%
}}}}
\put(3344,325){\makebox(0,0)[lb]{\smash{{\SetFigFont{11}{13.2}{\rmdefault}{\mddefault}{\updefault}{\color[rgb]{0,0,0}$m_v=0$}%
}}}}
\put(1711,-272){\makebox(0,0)[lb]{\smash{{\SetFigFont{11}{13.2}{\rmdefault}{\mddefault}{\updefault}{\color[rgb]{0,0,0}$\C^*.\tilde x$}%
}}}}
\put(3902,-153){\makebox(0,0)[lb]{\smash{{\SetFigFont{11}{13.2}{\rmdefault}{\mddefault}{\updefault}{\color[rgb]{0,0,0}$\O_X(-1)$}%
}}}}
\end{picture}%

%% file: slice4.pstex_t
\begin{picture}(0,0)%
\includegraphics{slice4.pstex}%
\end{picture}%
\setlength{\unitlength}{2901sp}%
\begingroup\makeatletter\ifx\SetFigFont\undefined%
\gdef\SetFigFont#1#2#3#4#5{%
  \reset@font\fontsize{#1}{#2pt}%
  \fontfamily{#3}\fontseries{#4}\fontshape{#5}%
  \selectfont}%
\fi\endgroup%
\begin{picture}(5855,3174)(1441,-3403)
\put(1441,-601){\makebox(0,0)[lb]{\smash{{\SetFigFont{11}{13.2}{\rmdefault}{\mddefault}{\updefault}{\color[rgb]{0,0,0}unstable orbit}%
}}}}
\put(6192,-476){\makebox(0,0)[lb]{\smash{{\SetFigFont{11}{13.2}{\rmdefault}{\mddefault}{\updefault}{\color[rgb]{0,0,0}$G$-orbit}%
}}}}
\put(6158,-986){\makebox(0,0)[lb]{\smash{{\SetFigFont{11}{13.2}{\rmdefault}{\mddefault}{\updefault}{\color[rgb]{0,0,0}$m^{-1}(0)$}%
}}}}
\put(6192,-1949){\makebox(0,0)[lb]{\smash{{\SetFigFont{11}{13.2}{\rmdefault}{\mddefault}{\updefault}{\color[rgb]{0,0,0}$K$-orbit}%
}}}}
\end{picture}%

%% file: interior.pstex_t
\begin{picture}(0,0)%
\includegraphics{interior.pstex}%
\end{picture}%
\setlength{\unitlength}{2368sp}%
\begingroup\makeatletter\ifx\SetFigFont\undefined%
\gdef\SetFigFont#1#2#3#4#5{%
  \reset@font\fontsize{#1}{#2pt}%
  \fontfamily{#3}\fontseries{#4}\fontshape{#5}%
  \selectfont}%
\fi\endgroup%
\begin{picture}(9508,4818)(2430,-4561)
\put(10951,-1636){\makebox(0,0)[lb]{\smash{{\SetFigFont{10}{12.0}{\rmdefault}{\mddefault}{\updefault}{\color[rgb]{0,0,0}$H_v$}%
}}}}
\put(3151,-4111){\makebox(0,0)[lb]{\smash{{\SetFigFont{10}{12.0}{\rmdefault}{\mddefault}{\updefault}{\color[rgb]{0,0,0}Stable}%
}}}}
\put(6001,-2536){\makebox(0,0)[lb]{\smash{{\SetFigFont{10}{12.0}{\rmdefault}{\mddefault}{\updefault}{\color[rgb]{0,0,0}Semistable}%
}}}}
\put(5251,-3736){\makebox(0,0)[lb]{\smash{{\SetFigFont{10}{12.0}{\rmdefault}{\mddefault}{\updefault}{\color[rgb]{0,0,0}$H_v$}%
}}}}
\put(5251, 89){\makebox(0,0)[lb]{\smash{{\SetFigFont{10}{12.0}{\rmdefault}{\mddefault}{\updefault}{\color[rgb]{0,0,0}$H_v$}%
}}}}
\put(2776, 14){\makebox(0,0)[lb]{\smash{{\SetFigFont{10}{12.0}{\rmdefault}{\mddefault}{\updefault}{\color[rgb]{0,0,0}= origin $\in\mathfrak t^*$}%
}}}}
\put(2776,-361){\makebox(0,0)[lb]{\smash{{\SetFigFont{10}{12.0}{\rmdefault}{\mddefault}{\updefault}{\color[rgb]{0,0,0}= weight $m\in\mathfrak t^*$}%
}}}}
\put(9151,-4561){\makebox(0,0)[lb]{\smash{{\SetFigFont{10}{12.0}{\rmdefault}{\mddefault}{\updefault}{\color[rgb]{0,0,0}Unstable}%
}}}}
\end{picture}%

%% file: smallfamily.pstex_t
\begin{picture}(0,0)%
\includegraphics{smallfamily.pstex}%
\end{picture}%
\setlength{\unitlength}{3315sp}%
\begingroup\makeatletter\ifx\SetFigFont\undefined%
\gdef\SetFigFont#1#2#3#4#5{%
  \reset@font\fontsize{#1}{#2pt}%
  \fontfamily{#3}\fontseries{#4}\fontshape{#5}%
  \selectfont}%
\fi\endgroup%
\begin{picture}(3575,2105)(901,-2293)
\put(2471,-368){\makebox(0,0)[lb]{\smash{{\SetFigFont{12}{14.4}{\rmdefault}{\mddefault}{\updefault}{\color[rgb]{0,0,0}$(\curly X_0,\L_0)$}%
}}}}
\put(3826,-368){\makebox(0,0)[lb]{\smash{{\SetFigFont{12}{14.4}{\rmdefault}{\mddefault}{\updefault}{\color[rgb]{0,0,0}$(\curly X_t,\L_t)\cong(X,L)$}%
}}}}
\put(4476,-747){\makebox(0,0)[lb]{\smash{{\SetFigFont{12}{14.4}{\rmdefault}{\mddefault}{\updefault}{\color[rgb]{0,0,0}$\forall t\ne0$}%
}}}}
\put(901,-376){\makebox(0,0)[lb]{\smash{{\SetFigFont{12}{14.4}{\rmdefault}{\mddefault}{\updefault}{\color[rgb]{0,0,0}$\curly X$}%
}}}}
\put(901,-2221){\makebox(0,0)[lb]{\smash{{\SetFigFont{12}{14.4}{\rmdefault}{\mddefault}{\updefault}{\color[rgb]{0,0,0}$\C$}%
}}}}
\end{picture}%

%% file: balanced.pstex_t
\begin{picture}(0,0)%
\includegraphics{balanced.pstex}%
\end{picture}%
\setlength{\unitlength}{1776sp}%
\begingroup\makeatletter\ifx\SetFigFont\undefined%
\gdef\SetFigFont#1#2#3#4#5{%
  \reset@font\fontsize{#1}{#2pt}%
  \fontfamily{#3}\fontseries{#4}\fontshape{#5}%
  \selectfont}%
\fi\endgroup%
\begin{picture}(7332,5349)(469,-5743)
\put(4606,-1336){\makebox(0,0)[lb]{\smash{{\SetFigFont{9}{10.8}{\rmdefault}{\mddefault}{\updefault}{\color[rgb]{0,0,0}$p_i$}%
}}}}
\put(7801,-1156){\makebox(0,0)[lb]{\smash{{\SetFigFont{11}{13.2}{\rmdefault}{\mddefault}{\updefault}{\color[rgb]{0,0,0}$\PP(H^0(L^r)^*)$}%
}}}}
\put(3676,-3886){\makebox(0,0)[lb]{\smash{{\SetFigFont{11}{13.2}{\rmdefault}{\mddefault}{\updefault}{\color[rgb]{0,0,0}$X$}%
}}}}
\put(5596,-1966){\makebox(0,0)[lb]{\smash{{\SetFigFont{9}{10.8}{\rmdefault}{\mddefault}{\updefault}{\color[rgb]{0,0,0}$\PP(\oplus_i\C_{p_i}^{m_i*})$}%
}}}}
\put(3001,-5536){\makebox(0,0)[lb]{\smash{{\SetFigFont{9}{10.8}{\rmdefault}{\mddefault}{\updefault}{\color[rgb]{0,0,0}$\PP(I^*)=\PP(H^0(L^r\otimes\I_{\cup_im_ip_i})^*)$}%
}}}}
\end{picture}%

%% file: dg1.pstex_t
\begin{picture}(0,0)%
\includegraphics{dg1.pstex}%
\end{picture}%
\setlength{\unitlength}{2901sp}%
\begingroup\makeatletter\ifx\SetFigFont\undefined%
\gdef\SetFigFont#1#2#3#4#5{%
  \reset@font\fontsize{#1}{#2pt}%
  \fontfamily{#3}\fontseries{#4}\fontshape{#5}%
  \selectfont}%
\fi\endgroup%
\begin{picture}(7152,2637)(946,-3094)
\put(1531,-961){\makebox(0,0)[lb]{\smash{{\SetFigFont{10}{12.0}{\rmdefault}{\mddefault}{\updefault}{\color[rgb]{0,0,0}$(X,L^r)$}%
}}}}
\put(3421,-601){\makebox(0,0)[lb]{\smash{{\SetFigFont{8}{9.6}{\rmdefault}{\mddefault}{\updefault}{\color[rgb]{0,0,0}Kempf-Ness}%
}}}}
\put(3691,-916){\makebox(0,0)[lb]{\smash{{\SetFigFont{8}{9.6}{\rmdefault}{\mddefault}{\updefault}{\color[rgb]{0,0,0}Zhang}%
}}}}
\put(4996,-781){\makebox(0,0)[lb]{\smash{{\SetFigFont{10}{12.0}{\rmdefault}{\mddefault}{\updefault}{\color[rgb]{0,0,0}Balanced $X\subset\PP^{N(r)}$}%
}}}}
\put(5716,-2671){\makebox(0,0)[lb]{\smash{{\SetFigFont{10}{12.0}{\rmdefault}{\mddefault}{\updefault}{\color[rgb]{0,0,0}cscK}%
}}}}
\put(6121,-1591){\makebox(0,0)[lb]{\smash{{\SetFigFont{8}{9.6}{\rmdefault}{\mddefault}{\updefault}{\color[rgb]{0,0,0}$r\to\infty$}%
}}}}
\put(1711,-2626){\makebox(0,0)[lb]{\smash{{\SetFigFont{10}{12.0}{\rmdefault}{\mddefault}{\updefault}{\color[rgb]{0,0,0}\textbf{???}}%
}}}}
\put(1576,-1861){\makebox(0,0)[lb]{\smash{{\SetFigFont{8}{9.6}{\rmdefault}{\mddefault}{\updefault}{\color[rgb]{0,0,0}HM\ \ criterion ?}%
}}}}
\put(946,-2806){\makebox(0,0)[lb]{\smash{{\SetFigFont{10}{12.0}{\rmdefault}{\mddefault}{\updefault}{\color[rgb]{0,0,0}For dim=0, multiplicity}%
}}}}
\put(991,-736){\makebox(0,0)[lb]{\smash{{\SetFigFont{10}{12.0}{\rmdefault}{\mddefault}{\updefault}{\color[rgb]{0,0,0}Stability of varieties}%
}}}}
\put(3961,-2536){\makebox(0,0)[lb]{\smash{{\SetFigFont{8}{9.6}{\rmdefault}{\mddefault}{\updefault}{\color[rgb]{0,0,0}???}%
}}}}
\put(7561,-2671){\makebox(0,0)[lb]{\smash{{\SetFigFont{10}{12.0}{\rmdefault}{\mddefault}{\updefault}{\color[rgb]{0,0,0}$Ham(X,\omega)$}%
}}}}
\put(7516,-781){\makebox(0,0)[lb]{\smash{{\SetFigFont{8}{9.6}{\rmdefault}{\mddefault}{\updefault}{\color[rgb]{0,0,0}$SU(N(r)\!+\!1)$}%
}}}}
\put(946,-3031){\makebox(0,0)[lb]{\smash{{\SetFigFont{10}{12.0}{\rmdefault}{\mddefault}{\updefault}{\color[rgb]{0,0,0}of any point $<\frac12$ total}%
}}}}
\put(3736,-1546){\makebox(0,0)[lb]{\smash{{\SetFigFont{8}{9.6}{\rmdefault}{\mddefault}{\updefault}{\color[rgb]{0,0,0}Donaldson}%
}}}}
\put(4906,-1771){\makebox(0,0)[lb]{\smash{{\SetFigFont{8}{9.6}{\rmdefault}{\mddefault}{\updefault}{\color[rgb]{0,0,0}Donaldson}%
}}}}
\put(3421,-1366){\makebox(0,0)[lb]{\smash{{\SetFigFont{8}{9.6}{\rmdefault}{\mddefault}{\updefault}{\color[rgb]{0,0,0}Tian}%
}}}}
\end{picture}%

%% file: dg2.pstex_t
\begin{picture}(0,0)%
\includegraphics{dg2.pstex}%
\end{picture}%
\setlength{\unitlength}{2901sp}%
\begingroup\makeatletter\ifx\SetFigFont\undefined%
\gdef\SetFigFont#1#2#3#4#5{%
  \reset@font\fontsize{#1}{#2pt}%
  \fontfamily{#3}\fontseries{#4}\fontshape{#5}%
  \selectfont}%
\fi\endgroup%
\begin{picture}(7501,2277)(991,-2464)
\put(4906,-511){\makebox(0,0)[lb]{\smash{{\SetFigFont{10}{12.0}{\rmdefault}{\mddefault}{\updefault}{\color[rgb]{0,0,0}Balanced $X\to Gr(N(r))$}%
}}}}
\put(7516,-511){\makebox(0,0)[lb]{\smash{{\SetFigFont{8}{9.6}{\rmdefault}{\mddefault}{\updefault}{\color[rgb]{0,0,0}$SU(N(r)\!+\!1)$}%
}}}}
\put(1036,-466){\makebox(0,0)[lb]{\smash{{\SetFigFont{10}{12.0}{\rmdefault}{\mddefault}{\updefault}{\color[rgb]{0,0,0}Stability of bundles}%
}}}}
\put(3421,-331){\makebox(0,0)[lb]{\smash{{\SetFigFont{8}{9.6}{\rmdefault}{\mddefault}{\updefault}{\color[rgb]{0,0,0}Kempf-Ness}%
}}}}
\put(3691,-646){\makebox(0,0)[lb]{\smash{{\SetFigFont{8}{9.6}{\rmdefault}{\mddefault}{\updefault}{\color[rgb]{0,0,0}Wang}%
}}}}
\put(6121,-1321){\makebox(0,0)[lb]{\smash{{\SetFigFont{8}{9.6}{\rmdefault}{\mddefault}{\updefault}{\color[rgb]{0,0,0}$r\to\infty$}%
}}}}
\put(1261,-2401){\makebox(0,0)[lb]{\smash{{\SetFigFont{10}{12.0}{\rmdefault}{\mddefault}{\updefault}{\color[rgb]{0,0,0}Slope criterion}%
}}}}
\put(5716,-2401){\makebox(0,0)[lb]{\smash{{\SetFigFont{10}{12.0}{\rmdefault}{\mddefault}{\updefault}{\color[rgb]{0,0,0}HYM}%
}}}}
\put(1576,-1501){\makebox(0,0)[lb]{\smash{{\SetFigFont{8}{9.6}{\rmdefault}{\mddefault}{\updefault}{\color[rgb]{0,0,0}HM\ \ criterion}%
}}}}
\put(4906,-1186){\makebox(0,0)[lb]{\smash{{\SetFigFont{8}{9.6}{\rmdefault}{\mddefault}{\updefault}{\color[rgb]{0,0,0}Donaldson}%
}}}}
\put(5176,-1411){\makebox(0,0)[lb]{\smash{{\SetFigFont{8}{9.6}{\rmdefault}{\mddefault}{\updefault}{\color[rgb]{0,0,0}Wang}%
}}}}
\put(1126,-1681){\makebox(0,0)[lb]{\smash{{\SetFigFont{8}{9.6}{\rmdefault}{\mddefault}{\updefault}{\color[rgb]{0,0,0}Mumford\ \,Gieseker}%
}}}}
\put(1576,-691){\makebox(0,0)[lb]{\smash{{\SetFigFont{10}{12.0}{\rmdefault}{\mddefault}{\updefault}{\color[rgb]{0,0,0}$E\to(X,L^r)$}%
}}}}
\put(3241,-2041){\makebox(0,0)[lb]{\smash{{\SetFigFont{8}{9.6}{\rmdefault}{\mddefault}{\updefault}{\color[rgb]{0,0,0}Donaldson-}%
}}}}
\put(3241,-2266){\makebox(0,0)[lb]{\smash{{\SetFigFont{8}{9.6}{\rmdefault}{\mddefault}{\updefault}{\color[rgb]{0,0,0}Uhlenbeck-Yau}%
}}}}
\put(7831,-2401){\makebox(0,0)[lb]{\smash{{\SetFigFont{10}{12.0}{\rmdefault}{\mddefault}{\updefault}{\color[rgb]{0,0,0}$U(E)$}%
}}}}
\put(991,-1861){\makebox(0,0)[lb]{\smash{{\SetFigFont{8}{9.6}{\rmdefault}{\mddefault}{\updefault}{\color[rgb]{0,0,0}Maruyama\ \ Simpson}%
}}}}
\end{picture}%

%% file: smallsubscheme.pstex_t
\begin{picture}(0,0)%
\includegraphics{smallsubscheme.pstex}%
\end{picture}%
\setlength{\unitlength}{3315sp}%
\begingroup\makeatletter\ifx\SetFigFont\undefined%
\gdef\SetFigFont#1#2#3#4#5{%
  \reset@font\fontsize{#1}{#2pt}%
  \fontfamily{#3}\fontseries{#4}\fontshape{#5}%
  \selectfont}%
\fi\endgroup%
\begin{picture}(5085,2592)(631,-3276)
\put(2296,-1456){\makebox(0,0)[lb]{\smash{{\SetFigFont{10}{12.0}{\rmdefault}{\mddefault}{\updefault}{\color[rgb]{0,0,0}$Z_1$}%
}}}}
\put(2071,-1951){\makebox(0,0)[lb]{\smash{{\SetFigFont{10}{12.0}{\rmdefault}{\mddefault}{\updefault}{\color[rgb]{0,0,0}$Z_2$}%
}}}}
\put(2206,-2266){\makebox(0,0)[lb]{\smash{{\SetFigFont{10}{12.0}{\rmdefault}{\mddefault}{\updefault}{\color[rgb]{0,0,0}$Z_0$}%
}}}}
\put(5716,-1366){\makebox(0,0)[lb]{\smash{{\SetFigFont{10}{12.0}{\rmdefault}{\mddefault}{\updefault}{\color[rgb]{0,0,0}$Z_1\times\Spec\C[t]/(t^2)$}%
}}}}
\put(5716,-1816){\makebox(0,0)[lb]{\smash{{\SetFigFont{10}{12.0}{\rmdefault}{\mddefault}{\updefault}{\color[rgb]{0,0,0}$Z_2\times\Spec\C[t]/(t^3)$}%
}}}}
\put(5716,-2401){\makebox(0,0)[lb]{\smash{{\SetFigFont{10}{12.0}{\rmdefault}{\mddefault}{\updefault}{\color[rgb]{0,0,0}$Z_0\times\{0\}$}%
}}}}
\put(4006,-1411){\makebox(0,0)[lb]{\smash{{\SetFigFont{10}{12.0}{\rmdefault}{\mddefault}{\updefault}{\color[rgb]{0,0,0}$X_0$}%
}}}}
\put(1801,-916){\makebox(0,0)[lb]{\smash{{\SetFigFont{12}{14.4}{\rmdefault}{\mddefault}{\updefault}{\color[rgb]{0,0,0}$X\times\C$}%
}}}}
\end{picture}%

%% file: smallnormalcone0.pstex_t
\begin{picture}(0,0)%
\includegraphics{smallnormalcone0.pstex}%
\end{picture}%
\setlength{\unitlength}{3315sp}%
\begingroup\makeatletter\ifx\SetFigFont\undefined%
\gdef\SetFigFont#1#2#3#4#5{%
  \reset@font\fontsize{#1}{#2pt}%
  \fontfamily{#3}\fontseries{#4}\fontshape{#5}%
  \selectfont}%
\fi\endgroup%
\begin{picture}(3577,3220)(1883,-3575)
\put(2693,-916){\makebox(0,0)[lb]{\smash{{\SetFigFont{9}{10.8}{\rmdefault}{\mddefault}{\updefault}{\color[rgb]{0,0,0}$X_t$}%
}}}}
\put(3829,-487){\makebox(0,0)[lb]{\smash{{\SetFigFont{9}{10.8}{\rmdefault}{\mddefault}{\updefault}{\color[rgb]{0,0,0}$Z_0'$}%
}}}}
\put(4037,-1815){\makebox(0,0)[lb]{\smash{{\SetFigFont{7}{8.4}{\rmdefault}{\mddefault}{\updefault}{\color[rgb]{0,0,0}$e$}%
}}}}
\put(1883,-1605){\makebox(0,0)[lb]{\smash{{\SetFigFont{9}{10.8}{\rmdefault}{\mddefault}{\updefault}{\color[rgb]{0,0,0}$Z_0\times\{t\}$}%
}}}}
\put(3380,-3526){\makebox(0,0)[lb]{\smash{{\SetFigFont{9}{10.8}{\rmdefault}{\mddefault}{\updefault}{\color[rgb]{0,0,0}$(\X)\_0=\widehat X\cup_eP$}%
}}}}
\end{picture}%

%% file: normalcone.pstex_t
\begin{picture}(0,0)%
\includegraphics{normalcone.pstex}%
\end{picture}%
\setlength{\unitlength}{4144sp}%
\begingroup\makeatletter\ifx\SetFigFont\undefined%
\gdef\SetFigFont#1#2#3#4#5{%
  \reset@font\fontsize{#1}{#2pt}%
  \fontfamily{#3}\fontseries{#4}\fontshape{#5}%
  \selectfont}%
\fi\endgroup%
\begin{picture}(3577,3254)(1883,-3575)
\put(2693,-916){\makebox(0,0)[lb]{\smash{\SetFigFont{12}{14.4}{\rmdefault}{\mddefault}{\updefault}{\color[rgb]{0,0,0}$X_t$}%
}}}
\put(4415,-540){\makebox(0,0)[lb]{\smash{\SetFigFont{12}{14.4}{\rmdefault}{\mddefault}{\updefault}{\color[rgb]{0,0,0}$Z_1'$}%
}}}
\put(3829,-487){\makebox(0,0)[lb]{\smash{\SetFigFont{12}{14.4}{\rmdefault}{\mddefault}{\updefault}{\color[rgb]{0,0,0}$Z_0'$}%
}}}
\put(4037,-1815){\makebox(0,0)[lb]{\smash{\SetFigFont{9}{10.8}{\rmdefault}{\mddefault}{\updefault}{\color[rgb]{0,0,0}$e$}%
}}}
\put(1959,-2105){\makebox(0,0)[lb]{\smash{\SetFigFont{12}{14.4}{\rmdefault}{\mddefault}{\updefault}{\color[rgb]{0,0,0}$Z_1\times\{t\}$}%
}}}
\put(1883,-1605){\makebox(0,0)[lb]{\smash{\SetFigFont{12}{14.4}{\rmdefault}{\mddefault}{\updefault}{\color[rgb]{0,0,0}$Z_0\times\{t\}$}%
}}}
\put(3380,-3526){\makebox(0,0)[lb]{\smash{\SetFigFont{12}{14.4}{\rmdefault}{\mddefault}{\updefault}{\color[rgb]{0,0,0}$(\X)\_0=\widehat X\cup_eP$}%
}}}
\end{picture}